\documentclass[10pt]{amsart}
\usepackage{amssymb}

\usepackage{amscd}
\usepackage{amsthm}
\usepackage[dvips]{graphicx}

\usepackage[all]{xy}

\newtheorem{Thm}{Theorem}[section]
\newtheorem{Lem}[Thm]{Lemma}

\newtheorem{Cor}[Thm]{Corollary}
\newtheorem{Prop}[Thm]{Proposition}

\newtheorem{Rem}[Thm]{Remark}

\title[On stable equivalences]{Classifying tame blocks and related algebras
up to stable equivalences of Morita type}
\author{ Guodong Zhou and Alexander Zimmermann}

\address{Guodong Zhou
\newline Institut f\"ur Mathematik,
\newline Universit\"at Paderborn,
\newline Warburger Stra\ss e 100,
\newline 33098 Paderborn,
\newline Germany}
\email{gzhou@math.uni-paderborn.de}
\address{Alexander Zimmermann
\newline Universit\'e de Picardie,
\newline D\'epartement de Math\'ematiques et LAMFA (UMR 6140 du CNRS),
\newline 33 rue St Leu,
\newline F-80039 Amiens Cedex 1,
\newline France}
\email{alexander.zimmermann@u-picardie.fr}
\date{March 17, 2010}

\newenvironment{Proof}[1][Proof]{\begin{trivlist}
\item[\hskip \labelsep {\bfseries #1}]}{\flushright
$\Box$\end{trivlist}}

\newcommand{\Char}{\mbox{char}}

\newcommand{\lra}{\longrightarrow}

\newcommand{\ra}{\rightarrow}
\newcommand{\sdp}{\times\kern-.2em\vrule height1.1ex depth-.05ex}
\newcommand{\epi}{\lra \kern-.8em\ra}

\newcommand{\N}{{\mathbb N}}

\newcommand{\ul}{\underline}

\newcommand{\Z}{{\mathbb Z}}

 \normalbaselines
\begin{document}

\begin{abstract}
We contribute to the classification of finite dimensional algebras
under stable equivalence of Morita type. More precisely
we give a classification
of the class of Erdmann's algebras of dihedral, semi-dihedral and
quaternion type and obtain as byproduct the validity of the
Auslander-Reiten conjecture for these classes of algebras.
\end{abstract}

\maketitle

\section*{Introduction}

Stable categories were introduced very early in the representation theory
of algebras and played a major r\^ole in the development of
Auslander-Reiten theory for example. Nevertheless, already in the 1970's
Auslander and Reiten knew that equivalences of stable categories can behave
very badly. For example there are indecomposable
finite dimensional algebras which are stably equivalent to a direct product
of two algebras none of which is separable \cite[Example 3.5]{AR}.

Around 1990 the concept of derived categories became popular in the representation
theory of groups and algebras by mainly two developments. First Happel interpreted
successfully tilting theory in the framework of derived categories and secondly
Brou\'e formulated his famous abelian defect group conjecture in this framework.
Many homological constructions are more natural in the language of derived
categories.
Work of Rickard \cite{Rickard1989} and Keller-Vossieck \cite{KellerVossieck}
show that an equivalence between derived
categories of self-injective algebras imply an equivalence between the stable
categories of these algebras of a very particular shape. They are induced by
bimodules which are invertible almost as for Morita equivalences. This discovery
in mind Brou\'e defined two algebras $A$ and $B$ to be stably equivalent of
Morita type if there is an $A-B$-bimodule $M$ and a $B-A$-bimodule $N$, which
are projective considered as module on either side only and so that there are
isomorphisms of bimodules
$M\otimes_BN\simeq A\oplus P$ for a projective $A-A$-bimodule $P$
and $N\otimes_AM\simeq B\oplus Q$ for a projective $B-B$-bimodule $Q$.

It soon became clear that stable equivalences of Morita type are much better
behaved than abstract stable equivalences. Nevertheless, classes of algebras
which are classified up to stable equivalence of Morita type are rare.
In recent joint work with Yuming Liu \cite{LZZ} we gave several invariants which
proved to be very sophisticated and powerful so that a classification of big
classes of symmetric algebras up to stable equivalence of Morita type becomes
feasible. The additional problem mainly is that the number of simple modules
is not proven to be an invariant under stable equivalence of Morita type. This
fact is the long-standing open Auslander-Reiten conjecture.

Erdmann gave an (up to parameters) finite list
\cite{ErdmannLNM} of algebras which are defined by properties on
their Auslander Reiten quiver and which include all blocks of
finite groups of tame representation type. Her classification is
up to Morita equivalence. Holm pursued further this approach and
classified the algebras in Erdmann's list up to derived
equivalence \cite{Holmhabil}. We shall give an account of his
results in Section~\ref{holmresults}.

In the present work we classify the algebras of dihedral, semi-dihedral and
quaternion type up to stable equivalence of Morita type. Our classification
is almost as complete as for derived equivalences and the classification coincides
in some sense with the derived equivalence classification. In particular we show the
Auslander-Reiten conjecture for these classes of algebras and note that
the classes are closed under stable equivalences of Morita type.

The paper is organised as follows. In Section~\ref{stableinv} we recall
some of the invariants under stable equivalence of Morita type we use in the
sequel. Section~\ref{holmresults} explains Holm's derived equivalence classification
of algebras of dihedral, semi-dihedral and quaternion type. In Section~\ref{Sectiontameblocks}
we present an independent classification for the case of tame blocks of group rings.
The proof is much simpler than the general case, and hence we decided to present the
arguments separately, though, of course, the general theorem includes this case as well.
Moreover, a short summary of Holm's result on Hochschild cohomology of tame blocks is given
there.
Section~\ref{dihedraltypesection}
shows that the derived equivalence classification of dihedral type algebras
coincides with the classification up to stable equivalence of Morita type.
The main tool is a result of Pogorza\l y \cite[Theorem 7.3]{Pogorzaly}.
This section is the first technical
core of the paper. Section~\ref{centreSection}
computes the centres of the algebras of semi-dihedral and of quaternion type.
This section prepares the classification result for these classes of algebras.
Section~\ref{thecentralpart} distinguishes then stable equivalence classes
of Morita type using basically invariants derived from the centre.
This part is the second technical core of the paper. Section~\ref{thestatement}
finally summarises large parts of what was proved before and contains
the main result Theorem~\ref{main}
of the paper as well as some results derived from K\"ulshammer
like invariants computed initially to distinguish derived equivalence classes.

\section{Stable invariants}
\label{stableinv}

The stable category $A-\ul{mod}$ of a finite dimensional $K$-algebra
$A$ has the same objects as the category
of $A$-modules and morphisms, denoted by $\ul{Hom}_A(M,N)$ from $M$ to $N$,
are equivalence classes of morphisms of $A$-modules modulo those factoring through a
projective $A$-module.

In this section we shall explain and state most of
the various properties of algebras
invariant under stable equivalences of Morita type used in the sequel.

The first reduction is a result of Keller-Vossieck and Rickard.

\begin{Thm} \label{KVR}(Keller-Vossieck~\cite{KellerVossieck} and
Rickard~\cite{Rickard1989})
Let $K$ be a field and let $A$ and $B$ be two self-injective $K$ algebras.
If the bounded derived categories of $A$ and $B$ are equivalent,
$D^b(A)\simeq D^b(B)$, then the algebras $A$ and $B$ are stably
equivalent of Morita type.
\end{Thm}

Hence in order to give a classification of a class of algebras
up to stable equivalence of Morita type we can start from a
classification up to derived equivalence and decide for two
representatives of the derived equivalence classes whether they are
stably equivalent of Morita type.

In order to do so we use several criteria, some
linked to questions around the centre
of the algebras.

We first recall a construction due to Brou\'e.
Let $A$ and $B$ be $K$-algebras. If $A$ is stably equivalent of Morita type
to $B$, then the subcategory of the stable category generated by
left and right projective $A\otimes_KA^{op}$-modules
is equivalent to the analogous category of
$B\otimes B^{op}$-modules. The $A\otimes_KA^{op}$-module $A$ is mapped to the
$B\otimes_KB^{op}$-module $B$ under this equivalence.
Therefore
$$\ul{End}_{A\otimes_KA^{op}}(A)\simeq \ul{End}_{B\otimes_KB^{op}}(B).$$
Brou\'e denotes by $$Z^{st}(A)=\ul{End}_{A\otimes_KA^{op}}(A)$$
the stable centre and by
$$Z^{pr}(A):=
\ker\left({End}_{A\otimes_KA^{op}}(A)
\lra\ul{End}_{A\otimes_KA^{op}}(A)\right)$$
the projective centre of $A$.

\begin{Thm} (Brou\'e \cite[Proposition
5.4]{Broue1994})\label{Broue} Let $A$ and $B$ be two algebras
which are stably equivalent of Morita type, then $Z^{st}(A)\simeq
Z^{st}(B)$.
\end{Thm}

The centre is usually not an invariant under stable equivalences of Morita type.
However one of the main results of \cite{LZZ} gives a partial answer.

\begin{Thm} (Liu, Zhou, Zimmermann \cite[Theorem 1.1]{LZZ})
Let $K$ be an algebraically closed field and let
$A$ and $B$ be two indecomposable
finite dimensional $K$-algebras which are stably equivalent
of Morita type. Then $\dim_K(HH_0(A))=\dim_K(HH_0(B))$ if and only if
the number of simple $A$-modules up to isomorphism
equals the number of simple $B$-modules up to isomorphism.
\end{Thm}

Since for symmetric algebras $Hom_K(HH_0(A),K)\simeq Z(A)$, we get that
for symmetric indecomposable finite dimensional algebras over
algebraically closed fields the dimension of the centres coincide
if and only if the number of simple modules coincide.

Moreover, a very useful criterion was given in \cite{LZZ} as well in order to
estimate the dimension of the projective centre.

\begin{Prop}\label{pRankCartan}
(Liu, Zhou, Zimmermann \cite[Proposition 2.4 and Corollary 2.9]{LZZ})
Let $K$ be an algebraically closed field and let $A$ be an
indecomposable symmetric $K$-algebra with $n$ simple modules up to isomorphism.
Then the dimension of the projective centre
equals the rank of the Cartan matrix, seen as linear mapping $K^n\lra K^n$.
\end{Prop}

A classical invariant, popularised recently by K\"ulshammer
\cite{Kuelshammerenglish}, is the Reynolds ideal defined for any
$K$-algebra as $R(A):=Z(A)\cap soc(A)$. For symmetric algebras $A$
and a perfect field $K$ of strictly positive characteristic
K\"ulshammer constructed a descending sequence of ideals
$T_n(A)^\perp$ of the centre of $A$, for $n\in\N$ with
$R(A)=\bigcap_{n\in\N}T_n^\perp(A)$.

\begin{Prop} \label{Higman}
\cite[Proposition 2.4 and proof of Proposition 2.5]{LZZ}
The projective centre of an algebra equals the Higman ideal of $A$ and
the Higman ideal of an algebra is in the socle of the algebra.
\end{Prop}

We shall use the following fact.

\begin{Thm} \label{Reynolds}
 Let $K$ be an algebraically closed field and let $A$ and $B$
be two finite dimensional symmetric indecomposable $K$-algebras
which are     stably equivalent of Morita type.     Then
 $$\dim_K(Z(A)/R(A))=\dim_K(Z(B)/R(B)).$$ Furthermore if $K$ is of positive
 characteristic or  the Cartan matrix of $A$ is non singular,
 then we have an isomorphism of
 algebras   $Z(A)/R(A)\simeq  Z(B)/R(B)$.
\end{Thm}

\begin{Proof} The first statement is the dual of  \cite[Corollary
5.4]{LZZ}, as all algebras in question are symmetric. For the
second statement, the case  of positive
 characteristic is contained in  \cite[Proposition 5.8]{KLZ}. In
 case of non singular Cartan matrix, by \cite[Proposition 5.1]{Xi2008} (or see the discussions at the end of this section), $B$ has
 also non  singular Cartan matrix. So the rank is the Cartan
 matrix  is equal to the number of
simple modules, that is,  the dimension of the Reynolds ideal. So
we have $Z^{pr}(A)=R(A)$ and $Z^{st}(A)=Z(A)/R(A)$. Now use
Theorem~\ref{Broue}.

\end{Proof}

Let $A$ be an indecomposable finite dimensional algebra and let $C_A$
be its Cartan matrix. The Cartan matrix induces in a natural way
a mapping of the Grothendieck group  $G_0(A)$ of abelian groups (the
Grothendieck group taken in the sense of $A$-modules modulo exact sequences).
The stable Grothendieck group $G_0^{st}(A)$ is defined as the cokernel of
$$G_0(A)\stackrel{ C_A}{\lra}G_0(A)\lra G_0^{st}(A)\lra 0$$

\begin{Prop} (Xi \cite{Xi2008})
Let $A$ and $B$ be finite dimensional indecomposable $K$-algebras and
suppose that $A$ and $B$ are stably equivalent of Morita type. Then
$G_0^{st}(A)\simeq G_0^{st}(A)$.
\end{Prop}

It is clear by this statement that a stable equivalence of Morita type
preserves those elementary divisors of the Cartan matrix which are different
from $1$, including their multiplicity. (Note that as usual the elementary
divisors are supposed to be non negative). In particular the absolute value of the
Cartan determinant is preserved.

\section{Algebras of dihedral, semi-dihedral and quaternion type}
\label{holmresults}

Let $K$ be an algebraically closed field. In this section we shall
give Karin Erdmann's list of algebras of dihedral, semi-dihedral
and quaternion type.

By Theorem~\ref{KVR} of Keller-Vossieck and Rickard,
for two self-injective algebras
$A$ and $B$, an equivalence $D^b(A)\simeq D^b(B)$ of the bounded derived
categories implies that $A$ and $B$ are stably equivalent of
Morita type. Hence, as basis of our discussion we shall use the
list of Thorsten Holm \cite{Holmhabil} of
algebras of dihedral, semi-dihedral and quaternion type up to
derived equivalences. There are three families: the algebras of
dihedral type, the algebras of semi-dihedral type, the algebras of
quaternion type. Each family is subdivided into three
subclasses: algebras with one simple module, algebras with two
simple modules and algebras with three simple modules. Each
subfamily contains
algebras with quivers and relations, depending on parameters.
$$\begin{array}{c||c|c|c}
&\mbox{dihedral}&\mbox{semidihedral}&\mbox{quaternion}\\ \hline
\mbox{1 simple}&K[X,Y]/(XY,X^m-Y^n),&SD(1\mathcal{A})_1^k,
k\geq 2;&Q(1\mathcal{A})_1^k, k\geq 2;\\
& m\geq n\geq 2, m+n>4;&&\\
&&&\\
&D(1\mathcal{A})_1^{1}=K[X,Y]/(X^2,Y^2); &&  \\
&&&\\
& (char K=2)   &(char(K)=2)\ SD(1\mathcal{A})_2^k(c,d)
&(char K=2)\ Q(1\mathcal{A})_2^k(c,d),\\
&K[X,Y]/(X^2,YX-Y^2);&k\geq 2, (c,d)\neq (0,0);&k\geq 2, (c,d)\neq (0,0);\\
&&&\\
&D(1\mathcal{A})^k_1, k\geq 2;&  & \\
&&&\\
&&&\\
&(char K=2)\ D(1\mathcal{A})^k_2(d),&&\\
&k\geq 2, d=0 \ or\  1;& &\\ \hline
\mbox{2 simples}&D(2\mathcal{B})^{k,s}(c),
&SD(2\mathcal{B})^{k,t}_1(c)&Q(2\mathcal{B})_1^{k,s}(a,c)\\
&k\geq s\geq 1, c\in\{0,\  1\}&k\geq 1, t\geq 2, c\in\{0, \ 1\};
&k\geq 1, s\geq 3, a\neq 0; \\
&&&\\
&&SD(2\mathcal{B})^{k,t}_2(c)&\\
& &k\geq 1, t\geq 2,  &\\
& &  k+t\geq 4,  c\in\{0, \ 1\};&\\
\hline
\mbox{3 simples}&D(3\mathcal{K})^{a,b,c},
&SD(3\mathcal{K})^{a,b,c}&Q(3\mathcal{K})^{a,b,c}\\
&a\geq b\geq c\geq 1;& a\geq b\geq c\geq 1, a\geq 2;
&a\geq b\geq c\geq 1, b\geq 2,\\
&&& (a, b, c)\neq (2,2,1);\\
&&&\\
 &D(3\mathcal{R})^{k,s,t,u},&&Q(3\mathcal{A})^{2,2}_1(d)\\
 &s\geq t\geq u\geq k\geq 1, t\geq 2&&d\not\in\{ 0, 1\}\\
\hline
\end{array}$$
All algebras with one simple module in the above list has the
quiver of type $1\mathcal{A}$
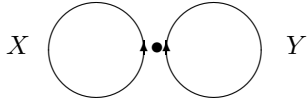
$$\unitlength0.6cm
\begin{picture}(7,3)
 \put(5.4,2){\circle{2.0}}
\put(2.8,2){\circle{2.0}} \put(3.85,1.9){\vector(0,1){0.3}}
\put(4,1.9){$\bullet$} \put(4.35,1.9){\vector(0,1){0.3}}
\put(0.8,1.9){$X$} \put(7,1.9){$Y$}
\end{picture}$$
and with relations
\begin{eqnarray*}
D(1\mathcal{A})^k_1 &: &  X^2, Y^2,(XY)^k-(YX)^k;\\
D(1\mathcal{A})^k_2(d)&: &  X^2-(XY)^k,Y^2-d\cdot(XY)^k,(XY)^k-(YX)^k,(XY)^kX,(YX)^kY;\\
SD(1\mathcal{A})_1^k&: &(XY)^k-(YX)^k,(XY)^kX,Y^2,X^2-(YX)^{k-1}Y;\\
SD(1\mathcal{A})_2^k(c,d)&:&(XY)^k-(YX)^k,(XY)^kX,Y^2-d(XY)^k,\\
&&\phantom{K<X,Y>/}X^2-(YX)^{k-1}Y+c(XY)^k;\\
Q(1\mathcal{A})_1^k&:&(XY)^k-(YX)^k,(XY)^kX,Y^2-(XY)^{k-1}X,X^2-(YX)^{k-1}Y;\\
Q(1\mathcal{A})_2^k(c,d)&:&X^2-(YX)^{k-1}Y-c(XY)^k,Y^2-(XY)^{k-1}X-d(XY)^k,\\
&&\phantom{K<X,Y>/}(XY)^k-(YX)^k,(XY)^kX,(YX)^kY.
\end{eqnarray*}

 The quivers of the algebras of type $2\mathcal{B}$, $3\mathcal{K}$,
 $3\mathcal{A}$ and $3\mathcal{R}$ are respectively:

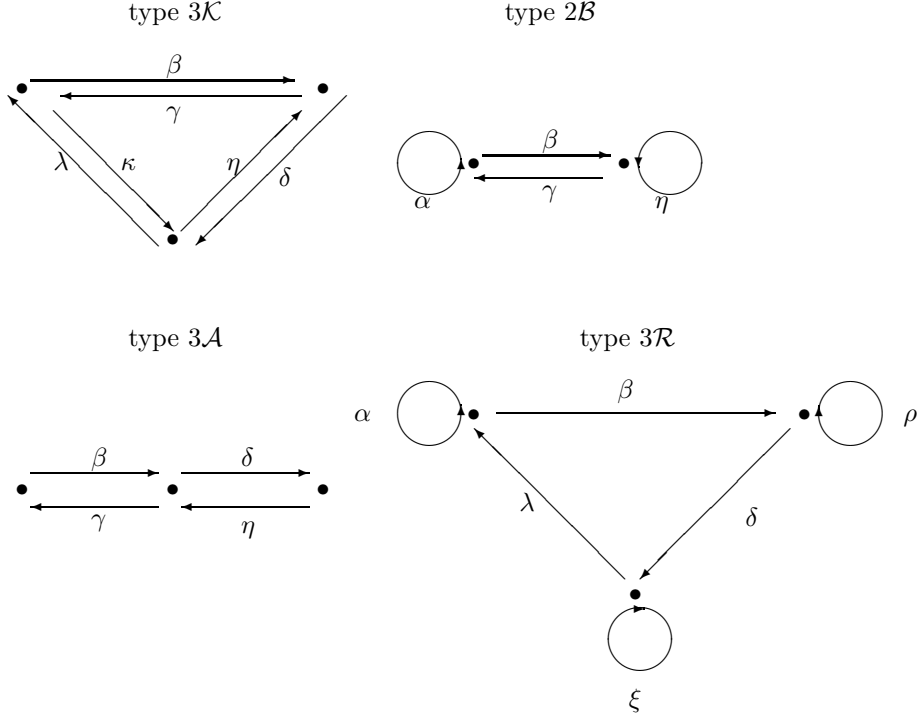
\unitlength1cm
\begin{picture}(15,5)
\put(3.5,4){type $3\mathcal{K}$} \put(2,3){$\bullet$}
\put(4,1){$\bullet$} \put(6,3){$\bullet$}
\put(2.2,3.2){\vector(1,0){3.5}}
\put(5.8,2.99){\vector(-1,0){3.2}}
\put(2.5,2.8){\vector(1,-1){1.6}} \put(3.9,1){\vector(-1,1){2}}
\put(4.2,1.2){\vector(1,1){1.6}} \put(6.4,3){\vector(-1,-1){2}}
\put(4,3.3){$\beta$} \put(4,2.7){$\gamma$} \put(5.5,1.8){$\delta$}
\put(4.8,2){$\eta$} \put(2.5,2){$\lambda$} \put(3.4,2){$\kappa$}

\put(8.5,4){type $2\mathcal{B}$} \put(8,2){$\bullet$}
\put(10,2){$\bullet$} \put(8.2,2.2){\vector(1,0){1.7}}
\put(9.8,1.9){\vector(-1,0){1.7}} \put(7.5,2.1){\circle{.8}}
\put(10.7,2.1){\circle{.8}} \put(7.92,2.15){\vector(0,1){.01}}
\put(10.28,2.05){\vector(0,-1){.01}} \put(7.3,1.5){$\alpha$}
\put(10.5,1.5){$\eta$} \put(9,2.3){$\beta$} \put(9,1.65){$\gamma$}


\end{picture}

\bigskip

\unitlength1cm
\begin{picture}(18,6)

\put(3.5,6){type $3\mathcal{A}$} \put(2,4){$\bullet$}
\put(4,4){$\bullet$} \put(6,4){$\bullet$}
\put(2.2,4.3){\vector(1,0){1.7}}
\put(3.9,3.85){\vector(-1,0){1.7}}
\put(4.2,4.3){\vector(1,0){1.7}}
\put(5.9,3.85){\vector(-1,0){1.7}}
\put(5,3.5){$\eta$} \put(3,4.4){$\beta$} \put(3,3.6){$\gamma$}
\put(5,4.4){$\delta$}

\put(9.5,6){type $3\mathcal{R}$}
\put(8,5){$\bullet$}
\put(10.15,2.6){$\bullet$}
\put(12.4,5){$\bullet$}
\put(12.68,5.15){\vector(0,1){.05}}
\put(8.4,5.1){\vector(1,0){3.7}}
 \put(10.1,2.9){\vector(-1,1){2}}
 \put(12.3,4.9){\vector(-1,-1){2}}
  \put(10,5.3){$\beta$}
\put(11.7,3.6){$\delta$} 
\put(8.7,3.8){$\lambda$} 

\put(7.5,5.1){\circle{.8}}\put(7.92,5.2){\vector(0,1){.01}}

\put(13.1,5.1){\circle{.8}}\put(12.7,53.2){\vector(0,1){.01}}

\put(10.3,2.1){\circle{.8}}\put(10.35,2.5){\vector(1,0){.01}}

\put(6.5,5){$\alpha$} \put(10.15,1.2){$\xi$} \put(13.8,5){$\rho$}

\end{picture}

The relations are respectively
\begin{eqnarray*}
D(2\mathcal{B})^{k, s}(c) &: &
\beta\eta, \eta\gamma, \gamma\beta, \alpha^2-c(\alpha\beta\gamma)^k,
(\alpha\beta\gamma)^k-(\beta\gamma\alpha)^k, \eta^s-(\gamma\alpha\beta)^k;\\
SD(2\mathcal{B})^{k, t}_1(c)&: & \gamma\beta, \eta\gamma,
\beta\eta,
\alpha^2-(\beta\gamma\alpha)^{k-1}\beta\gamma-c(\alpha\beta\gamma)^k,
\eta^t-(\gamma\alpha\beta)^k,
(\alpha\beta\gamma)^k-(\beta\gamma\alpha)^k;\\
SD(2\mathcal{B})^{k, t}_2(c)&:&
\beta\eta-(\alpha\beta\gamma)^{k-1}\alpha\beta,
\eta\gamma-(\gamma\alpha\beta)^{k-1}\gamma\alpha,
\gamma\beta-\eta^{t-1}, \alpha^2-c(\alpha\beta\gamma)^k,
\beta\eta^2, \eta^2\gamma;\\
Q(2\mathcal{B})^{k, s}_1(a,c)&:&\gamma\beta-\eta^{s-1},
\beta\eta-(\alpha\beta\gamma)^{k-1}\alpha\beta,
\eta\gamma-(\gamma\alpha\beta)^{k-1}\gamma\alpha,
 \\
&&\phantom{K<X,Y>/}\alpha^2-a(\beta\gamma\alpha)^{k-1}\beta\gamma-c
(\beta\gamma\alpha)^k, \alpha^2\beta, \gamma\alpha^2;\\
D(3\mathcal{K})^{a, b, c} &: & \beta\delta, \delta\lambda,
\lambda\beta, \gamma\kappa, \kappa\eta, \eta\gamma,
(\beta\gamma)^a-(\kappa\lambda)^b,
(\lambda\kappa)^b-(\eta\delta)^c, (\delta\eta)^c-(\gamma\beta)^a;\\
D(3\mathcal{R})^{k, s, t, u} &: & \alpha\beta, \beta\rho,
\rho\delta, \delta\xi, \xi\lambda, \lambda\alpha,
\alpha^s-(\beta\delta\lambda)^k, \rho^t-(\delta\gamma\beta)^k,
\xi^u-(\lambda\beta\delta)^k  ;\\
SD(3\mathcal{K})^{a, b, c} &: & \kappa\eta, \eta\gamma, \gamma\kappa,
\delta\gamma-(\gamma\alpha)^{a-1}\gamma,
\beta\delta-(\kappa\lambda)^{b-1}\kappa, \lambda\beta-(\eta\delta)^{c-1}\eta;\\
Q(3\mathcal{K})^{a, b, c} &: &
\beta\delta-(\kappa\lambda)^{a-1}\kappa,
\eta\gamma-(\lambda\kappa)^{a-1}\lambda,
\delta\lambda-(\gamma\beta)^{b-1}\gamma,
\kappa\eta-(\beta\gamma)^{b-1}\beta,
\lambda\beta-(\eta\delta)^{c-1}\eta, \\
&&\phantom{K<X,Y>/} \gamma\kappa-(\delta\eta)^{c-1}\delta,
\gamma\beta\delta, \delta\eta\gamma, \lambda\kappa\eta;\\
Q(3\mathcal{A})^{2,2}_1(d) &: & \beta\delta\eta-\beta\gamma\beta,
\delta\eta\gamma-\gamma\beta\gamma,
\eta\gamma\beta-d\eta\delta\eta,
\gamma\beta\delta-d\delta\eta\delta, \beta\delta\eta\delta,
\eta\gamma\beta\gamma.
\end{eqnarray*}

The following result suggests that we only need to consider
internally these three classes of algebras in order to classify
them up to stable equivalences of Morita type.

\begin{Prop} \label{onlyinternal}
If two indecomposable algebras $A$ and $B$ are stably equivalent
of Morita type and $A$ is of dihedral (resp. semi-dihedral,
quaternion) type, then so is $B$.
\end{Prop}

\begin{Proof}
These classes of algebras are defined in terms of
the nature of their Auslander-Reiten quiver. Roughly, an
algebra $A$ is of one of these types if
\begin{itemize}
\item $A$ is symmetric, indecomposable and tame;

\item the Cartan matrix of $A$ is non-singular.

\item the stable Auslander Reiten quiver of $A$ has the following components

\begin{tabular}{l||l|l|l}
&dihedral type& semidihedral type& quaternion type\\ \hline\hline &&&\\
tubes& rank 1 and 3& rank at most 3& rank at most 2\\ &&&\\
&at most two $3$-tubes& at most one $3$-tube&\\ \hline &&&\\
others& $\Z A_\infty^\infty/\Pi$&$\Z A_\infty^\infty$ and $\Z D_\infty$\\
\hline
\end{tabular}

\medskip

 for a certain group $\Pi$.

\end{itemize}
For more details we refer to \cite{ErdmannLNM}.

By a result of Yuming Liu (\cite[Corollary 2.4]{Liu2008})(resp.
Henning Krause (\cite[last corollary of the article]{Krause})), if two algebras are stably
equivalent of Morita type and one of them is symmetric (resp.
tame), so is the other. If two algebras are stably equivalent of
Morita type, they are stably equivalent and thus their stable
Auslander-Reiten quiver are isomorphic. By a result of Chang-Chang
Xi (\cite[Proposition 5.1]{Xi2008}),  if two algebras are stably
equivalent of Morita type, the absolute values of the determinant
of their Cartan matrices are the same and thus if the Cartan
matrix of one algebra is non-singular, so is that of the other.
Therefore, the defining properties are preserved by a stable
equivalence of Morita type between two indecomposable algebras.
\end{Proof}

\section{tame blocks}

\label{Sectiontameblocks}

\subsection{Derived classification}

The following is a classification of algebras up to derived
equivalence, as given by Holm \cite{Holmhabil},   which could
occur as blocks of group algebras. For some cases the question if
there is a block of a group with this derived equivalence type is
not clear yet. We include in this case the algebra as well. Now
let $K$ be an algebraically closed field of characteristic two.
Let $A$ be a  tame block of defect $n\geq 2$. Then $A$ is derived equivalent
to one of the following algebras.
$$\begin{array}{c||c|c|c}
&\mbox{dihedral}&\mbox{semidihedral}&\mbox{quaternion}\\  \hline&&&\\
\mbox{1 simple}& D(1\mathcal{A})^{2^{n-2}}_1, n\geq 2;
&SD(1\mathcal{A})_1^{2^{n-2}}, n\geq 4;&Q(1\mathcal{A})_1^{2^{n-2}}, n\geq 3;\\
 \hline &&&\\
\mbox{2 simples}&D(2\mathcal{B})^{1,2^{n-2}}(c),
&SD(2\mathcal{B})^{1, 2^{n-2}}_1(c), & Q(2\mathcal{B})_1^{2, 2^{n-2}}(a,c), \\
&c\in\{0 , 1\}, n\geq 3;&c\in\{0 , 1\}, n\geq 4;&
n\geq 3, a\in K^*, c\in K;\\
&&&\\
 &&SD(2\mathcal{B})^{2, 2^{n-2}}_2(c),  &\\
  && c\in\{0, 1\}, n\geq 4;&\\
  \hline &&&\\
\mbox{3 simples}&D(3\mathcal{K})^{2^{n-2},1, 1}, n\geq 2;
&SD(3\mathcal{K})^{2^{n-2}, 2, 1}, n\geq 4;&
Q(3\mathcal{K})^{2^{n-2}, 2, 2}, n\geq 3.\\
\hline
\end{array}$$

\subsection{Hochschild cohomology of tame blocks}

\label{HochschildSection}

If $A$ and $B$ are two algebras which are stably equivalent of
Morita type, then Xi shows \cite[Theorem 4.2]{Xi2008} that the
Hochschild cohomology groups $HH^m(A)$ and $HH^m(B)$ are
isomorphic for any $m\geq 1$.  Furthermore, in a recent paper of
the first author with Shengyong Pan (\cite{PanZhou}), we proved
that a stable equivalence of Morita type preserve the algebra
structure of the stable Hochschild cohomology, that is, the
Hochschild cohomology modulo the projective center.

For the sake of completeness we resume results of Holm
\cite{Holmhabil} which allow to distinguish a certain number of
pairs of algebras up to stable equivalence of Morita type,
although we could avoid using  these results in the sequel, mainly
because they only deal with blocks of group rings with one or
three simple modules.

\subsubsection{Dihedral type}

By \cite[Theorem 3.2.2]{Holmhabil} the Hochschild cohomology ring
of a block with dihedral defect group of order $2^n$ with $n\geq
2$ and one simple module has dimension
$dim(HH^{i}(B))=2^{n-2}+3+4i$.

By \cite[Theorem 3.2.8]{Holmhabil} the Hochschild cohomology ring
of a block with dihedral defect group of order $2^n$ with $n\geq
2$ and three simple modules has dimension $2^{n-2}+3$ in degree
$0$, and dimension $2^{n-2}+1$ in degree $1$. Further, for all
$i\geq 1$, $dim(HH^{3i-1}(B))=2^{n-2}-1+4i$ and
$dim(HH^{3i}(B))=dim(HH^{3i+1}(B))=2^{n-2}+1+4i$.

\subsubsection{Semi-dihedral type}
By \cite[Theorem 3.3.2]{Holmhabil} the Hochschild cohomology ring
of a block with semi-dihedral defect group of order $2^n$ with
$n\geq 4$ and one simple module has dimension $2^{n-2}+3$ in
degree $0$, dimension $2^{n-2}+6$ in degree $1$, dimension
$2^{n-2}+7$ in degree $2$, and dimension $2^{n-2}+8$ in degree
$3$. Further, $dim(HH^{i+4}(B))=dim(HH^i(B))+8$.

By \cite[Theorem 3.3.3]{Holmhabil} the Hochschild cohomology ring
of a block with semi-dihedral defect group of order $2^n$ with
$n\geq 4$ and three simple modules has dimension $2^{n-2}+4$ in
degrees $0$ and $3$, dimension $2^{n-2}+2$ in degrees $1$ and $2$,
and dimension $2^{n-2}+5$ in degree $4$. Further,
$dim(HH^{i+4}(B))=dim(HH^i(B))+2+x(i)$, where $x(i)$ is $0$ if $3$
divides $i$, and $x(i)=1$ else.

\subsubsection{Quaternion type}
By \cite[Theorem 3.4.2]{Holmhabil} a block with one simple modules
and quaternion defect group of order $2^n$ with $n\geq 3$ has
periodic Hochschild cohomology ring with period $4$ and dimension
$2^{n-2}+3$ in degrees congruent $0$ or $3$ mod $4$ and of
dimension $2^{n-2}+5$ in degrees congruent $1$ or $2$ mod $4$.

By \cite[Theorem 3.4.6]{Holmhabil} a block with three simple
modules and quaternion defect group of order $2^n$ with $n\geq 3$
has periodic Hochschild cohomology ring with period $4$ and
dimension $2^{n-2}+5$ in degrees congruent $0$ or $3$ mod $4$ and
of dimension $2^{n-2}+3$ in degrees congruent $1$ or $2$ mod $4$.

\subsection{Blocks of dihedral defect groups}

\begin{Prop}  Let $K$ be an algebraically closed field of
characteristic $2$ and let $A$ be a dihedral block of defect $n\geq 2$.
Then $A$ is stably  equivalent of Morita type  to one and exactly one
of the following algebras:
$D(1\mathcal{A})_1^{2^{n-2}}$; $D(2\mathcal{B})^{1,2^{n-2}}(c)$ (for $n\geq 3$)
with $c=0$ or $c=1$; $D(3\mathcal{K})^{2^{n-2}, 1, 1}$.

As a
consequence,  the derived classification coincides with the
classification up to stable equivalences of Morita type.
\end{Prop}

\begin{Rem}\rm
Before giving the proof, we remark that
for a dihedral block with two simple modules, we don't know
whether the case $c=1$ really occurs. All known examples have zero
as the value of this scalar. But this doesn't influence our
result, since $D(2\mathcal{B})^{k, s}(0)$ is NOT derived
equivalent to $D(2\mathcal{B})^{k, s}(1)$. There are several
proofs of this fact (cf \cite[Corollary 5.3]{KauerRoggenkamp}
\cite[Theorem 1.1]{HolmZimmermann}). One can
also use a result of Pogorza\l y(\cite[Theorem 7.3]{Pogorzaly}),
which says that an algebra stably equivalent to a self-injective
special biserial algebra which is not a Nakayama algebra is itself
a selfinjective special biserial algebra. Notice that
$D(2\mathcal{B})^{k, s}(0)$ is a symmetric special biserial
algebra, but $D(2\mathcal{B})^{k, s}(1)$ is not. As a consequence
the algebras $D(2\mathcal{B})^{k, s}(1)$ cannot be stably
equivalent to any algebra of the other classes.
\end{Rem}

\begin{Proof} Since Holm's result \cite{Holmhabil} implies that
any algebra of dihedral type is derived equivalent to one in the list we gave,
we just need to show that any two algebras in the list are not stably equivalent
of Morita type.

We prove that for different parameter $s\neq t$,
$D(2\mathcal{B})^{1, s}(1)$ is NOT stably equivalent of Morita
type to $D(2\mathcal{B})^{1, t}(1)$. To this end, one computes the
dimension of the stable centre, that is, the quotient of the
centre by the projective centre. By Proposition~\ref{pRankCartan},
for a symmetric algebra, the dimension of the projective centre is
the $p$-rank of the Cartan matrix, where $p$ is the characteristic
of the ground field, which is two for tame blocks. We have thus
that for $A=D(2\mathcal{B})^{1, 2^{n-2}}(1)$, $dim Z^{st}(A)=2^{n-2}+2$
for $n\geq 3$. Since $n\geq 3$ this dimension distinguishes
two algebras with different parameters in this class.
Another way to see this
is to use the absolute value of the determinant of the
Cartan matrix, which is invariant under stable equivalences of
Morita type, by a result of Chang-Chang Xi(\cite[Proposition
5.1]{Xi2008}). In fact,  the absolute value of the determinant of
the Cartan matrix of $D(2\mathcal{B})^{1, s}(1)$ is $4s$.

Now consider other classes of algebras. Pogorza\l y proved the
Auslander-Reiten conjecture for self-injective special biserial
algebras (\cite[Theorem 0.1]{Pogorzaly}), that is, if two
self-injective special biserial algebras are stably equivalent,
they have the same number of non projective simple modules. Thus
two indecomposable non-simple self-injective special biserial
algebras with different numbers of simple modules cannot be
stably equivalent. By \cite[Corollary 1.2]{LZZ}, we know that for
symmetric algebras, this is equivalent to say that their centre have
the same dimension. Now by computing the dimension of the centre,
we obtain easily that the number of simple modules and the defect
$n$ characterise equivalence classes under stable equivalences of
Morita type  of dihedral blocks which are special biserial. On can
also use the computations of Holm about Hochschild cohomology of
dihedral blocks resumed in Section~\ref{HochschildSection}
to distinguish dihedral blocks with one simple
module from those with three simple modules.
\end{Proof}

\subsection{Blocks with semi-dihedral defect groups}

\begin{Prop}  Let $K$ be an algebraically closed field of
characteristic $2$ and let
 $A$ be a semi-dihedral block of defect $n\geq 4$.
 Then $A$ is stably  equivalent of Morita type  to one
 of the following algebras:
$SD(1\mathcal{A})^{2^{n-2}}$ with $n\geq 4$; $SD(2\mathcal{B})^{1,
2^{n-2}}_1(c)$ with $n\geq 4$, $c\in\{0,1\}$;
$SD(2\mathcal{B})^{2, 2^{n-2}}_2(c)$ with $n\geq 4$,
$c\in\{0,1\}$; $SD(3\mathcal{K})^{2^{n-2}, 2, 1}, n\geq 4$.
\end{Prop}

\begin{Rem} \rm
\begin{enumerate}
\item The list of algebras occurring as blocks of group algebras
is taken from \cite{Holmhabil}.

\item In the above classification we have two problems still.
There is a scalar problem, that is, as in the case of derived
equivalence classification, we cannot determine whether for
different values of $c$, $SD(2\mathcal{B})^{1, 2^{n-2}}_1(0)$
(resp. $SD(2\mathcal{B})^{2, 2^{n-2}}_2(0)$) is not stably
equivalent of Morita type to $SD(2\mathcal{B})^{1, 2^{n-2}}_1(1)$
(resp. $SD(2\mathcal{B})^{2, 2^{n-2}}_2(1)$).

Moreover, we do not know whether the two algebras
$SD(2\mathcal{B})^{1, 2^{n-2}}_1(c_1)$
and $SD(2\mathcal{B})^{1, 2^{n-2}}_2(c_2)$ are stably equivalent of Morita type.
Therefore, up to these problems, the derived classification coincides with the
classification up to stable equivalences of Morita type.
\end{enumerate}
\end{Rem}

\begin{Proof} Since a derived equivalence between self-injective
algebras induces a stable equivalence of Morita
type, the statement of the proposition is true simply by the
derived equivalence classification of Thorsten Holm. We now prove
that the classification is complete up to the problems cited above.

By the result of Thorsten Holm on Hochschild cohomology of
semi-dihedral blocks, a semi-dihedral block with one simple module
can not be stably equivalent of Morita type to a semi-dihedral
block with three simple modules. The dimension of the stable
centre of $SD(1\mathcal{A})^{2^{n-2}}$ with $n\geq 4$ is
$2^{n-2}+3$, it is $2^{n-2}+2$ for $SD(2\mathcal{B})^{1,
2^{n-2}}_1(c)$ and is $2^{n-2}+4$ for $SD(2\mathcal{B})^{2,
2^{n-2}}_2(c)$, while for $SD(3\mathcal{K})^{2^{n-2}, 2, 1}$, it
is $2^{n-2}+3$. This invariant distinguishes semi-dihedral blocks
with two simple modules from those with one or three simple
modules and it also distinguishes $SD(2\mathcal{B})^{1,
2^{n-2}}_1(c)$ from $SD(2\mathcal{B})^{2, 2^{n-2}}_2(c)$.

\end{Proof}

\subsection{Blocks with quaternion defect groups}

\begin{Prop}  Let $K$ be an algebraically closed field of
characteristic $2$ and let
 $A$ be a   block with generalised quaternion defect groups of
 defect $n\geq 3$.
 Then $A$ is stably  equivalent of Morita type  to one
 of the following algebras:
$Q(1\mathcal{A})^{2^{n-2}} $ with $n\geq 3$;
$Q(2\mathcal{B})_1^{2, 2^{n-2}}(a,c)$ with $n\geq 3, a\in K^*,
c\in K$;   $Q(3K)^{2^{n-2}, 2, 2}$ with  $n\geq 3$.
\end{Prop}

\begin{Rem}\rm The above classification is complete up to some scalar
problem, that is, as in the case of derived equivalence
classification, we cannot determine whether for different values
of $a$ and $c$, $Q(2\mathcal{B})_1^{2, 2^{n-2}}(a,c)$   is not
stably equivalent of Morita type to $Q(2\mathcal{B})_1^{2,
2^{n-2}}(a',c')$. Therefore, up to these scalar problems,   the
derived classification coincides with the classification up to
stable equivalences of Morita type.
\end{Rem}

\begin{Proof} Since a derived equivalence between self-injective
algebras induces a stable equivalence of Morita
type, the statement of the proposition is true simply by the
derived equivalence classification of Thorsten Holm. We now prove
that the classification is complete up to the scalar problem.

The dimension of the stable
centre  is $2^{n-2}+3$ for $Q(1\mathcal{A})^{2^{n-2}} $,  is
$2^{n-2}+4$ for $Q(2\mathcal{B})_1^{2, 2^{n-2}}(a,c)$ and is
$2^{n-2}+5$ for $Q(3\mathcal{K})^{2^{n-2}, 2, 2}$. This invariant thus
distinguishes these algebras up to stable equivalences of Morita
type up to the scalar problem. One can also use the result of
Thorsten Holm on Hochschild cohomology of blocks with generalised
quaternion defect groups to distinguish blocks with generalised
quaternion defect groups having one simple module from those having
three simple modules.
\end{Proof}

\section{Algebras of dihedral type}
\label{dihedraltypesection}

We classify algebras of dihedral type up to stable equivalences of
Morita type in this section.  Notice that  all algebras except
$B_1=K[X, Y]/(X^2, Y^2-XY)$ and $D(1\mathcal{A})_2^k(d)$ are
special biserial. By the result of Pogorza\l y \cite[Theorem 0.1]{Pogorzaly},
one only needs to consider separately dihedral algebras with one, two
or three simple modules.

\subsection{One  simple module}
Let $K$ be an algebraically closed field of characteristic $p\geq
0$. By the classification of Erdmann (\cite{ErdmannLNM}), an
algebra of dihedral type with one simple module is Morita
equivalent to one of the following algebras.

\begin{center}
$A_1(m, n)=K[X,Y]/(XY,X^m-Y^n)$ with
$m\geq n\geq 2$ and $m+n>4$;
\end{center}

\begin{center}
$C_1=K[X, Y]/(X^2, Y^2]$;
\end{center}

\begin{center}
$D(1\mathcal{A})_1^k=K\langle X, Y\rangle/(X^2, Y^2, (XY)^k-(YX)^k$
with $k\geq 2$;
\end{center}

and if $p=2$,
\begin{center}
$B_1=K[X, Y]/(X^2, Y^2-XY)$;
\end{center} and

\begin{center}
$D(1\mathcal{A})_2^k(d)= K\langle X,Y\rangle/(X^2-(XY)^k,Y^2-d(XY)^k,
(XY)^k-(YX)^k,(XY)^kX,(YX)^kY)$.
\end{center}

\begin{Prop}\label{dihedralonesimple}
Let $K$ be an algebraically closed field of characteristic
$p\geq 0$ and let $A$ be an algebra
of dihedral type with one simple module. Then $A$ is stably
equivalent of Morita type to one and exactly one  algebra in the
following list:
\begin{itemize}
\item $A_1(n, m)$ with $m\geq n\geq 2$ and $m+n>4$;
\item $C_1 $;
\item $D(1\mathcal{A})_1^k $ with $k\geq 2$;
\item  if $p=2$, $B_1$  and $D(1\mathcal{A})_2^k(d)$
with $k\geq 2$ and $d\in \{0, 1\}$, except that we don't
  know whether $D(1A)^k_2(0) $ and  $D(1A)^k_2(1)$ are
  stably equivalent of Morita type or not.
\end{itemize}
\end{Prop}

The proof combines the following five claims below using
some invariants of these algebras shown in the following table.

Characteristic zero case
$$\begin{array}{c||c|c|c|c}
\mbox{algebra}\  A&  A_1(m, n)&C_1& D(1\mathcal{A})_1^k  \\ \hline
 dim Z(A) & n+m& 4& k+3\\
dim Z^{pr}(A) &  1 & 1& 1\\
dim Z^{st}(A)& n+m-1& 3 & k+2\\
C_A &[n+m] & [4] & [4k]\\
G_0^{st}&\mathbb{Z}/{(n+m)} &  \mathbb{Z}/4 & \mathbb{Z}/{4k}  \\
  \hline
\end{array}$$
Characteristic two case
$$\begin{array}{c||c|c|c|c|c|c}
\mbox{algebra } A&  A_1(m, n)&C_1& D(1{\mathcal A}_1^k) & B_1&
D(1\mathcal{A})_2^k(d)  \\ \hline
 dim Z(A) & n+m& 4& k+3 &4 & k+3\\
dim Z^{pr}(A) &  0\  \mbox{or}\  1 & 0& 0& 0 & 0 \\
dim Z^{st}(A)& n+m\ \mbox{or}\ n+m-1  & 4 & k+3& 4& k+3\\
C_A &[n+m] & [4] & [4k]& [4] & [4k]\\
G_0^{st}&\mathbb{Z}/{(n+m)} &  \mathbb{Z}/4 & \mathbb{Z}/{4k} &
\mathbb{Z}/4 & \mathbb{Z}/{4k}  \\
  \hline
\end{array}$$
Characteristic $p>2$ case
$$\begin{array}{c||c|c|c|}
\mbox{algebra } A&  A_1(m, n)&C_1& D(1\mathcal{A})_1^k   \\ \hline
 dim Z(A) & n+m& 4& k+3 \\
dim Z^{pr}(A) &  0\  \mbox{or}\  1 & 1& 0\  \mbox{or}\  1\\
dim Z^{st}(A)& n+m\mbox{ or }n+m-1& 3 & k+3\mbox{ or }k+2\\
C_A &[n+m] & [4] & [4k]\\
G_0^{st}&\mathbb{Z}/{(n+m)} &  \mathbb{Z}/4 & \mathbb{Z}/{4k}   \\
  \hline
\end{array}$$

By the result of Pogorza\l y (\cite[Theorem 7.3]{Pogorzaly}),
 we only need to compare $A_1(m, n), C_1$ with $D_1(k)$,
since they are special biserial, and compare $B_1$ with
$D_1(1A)^k(d)$, since they are not special biserial.

\bigskip

\textit{Claim 1. $C_1$  cannot be stably equivalent of Morita type
to $A_1(m, n)$ or  $D(1\mathcal{A})_1^k $.}

Comparing the stable Grothendieck groups gives the result. Since $m+n> 4$
and $k\geq 2$.

Similarly one proves

\bigskip

\textit{Claim 1'. $B_1$  cannot be stably equivalent of Morita
type to $D(1\mathcal{A})_2^k (d)$.}

\bigskip

\textit{Claim 2. $A_1(m, n)$ cannot be stably equivalent of Morita
type to $D(1\mathcal{A})_1^k $.  }

Compare  their stable centres    and  their stable Grothendieck
groups.

\bigskip

\textit{Claim 3.  $A=A_1(m, n)$ is not stably equivalent of Morita
type to $A'=A_1(m', n')$ for $(m, n)\neq (m', n')$ }

Now suppose that $A=A_1(m, n)$ is stably equivalent of Morita type
to $A'=A(m', n')$, then by comparing their stable Grothendieck
groups, $n+m=n'+m'$.
The Loewy length of the stable centre of $A(m, n)$ is $max(m,
n)$ if the characteristic $p$ divides $m+n$ and is $max(m, n)+1$, otherwise.  Thus
the stable centres are not isomorphic.

\bigskip

\textit{Claim 4.   $D(1\mathcal{A})_1^k $ cannot be stably
equivalent of Morita type to    $D(1\mathcal{A})_1^l $ for $k\neq
l$}

Comparing the orders of the stable Grothendieck groups gives the
result.

\bigskip

\textit{Claim 5. $D(1\mathcal{A})^k_2(d)$  cannot be stably
equivalent of Morita type to $D(1\mathcal{A})^l_2(d)$ for $k\neq
l$.}

 Consider the stable Grothendieck groups or the stable centres.

\bigskip

\subsection{Two  simple modules}

For algebras of dihedral type with two simple modules, we have the
following result of Holm.

\begin{Prop}(\cite[Proposition 2.3.1]{Holmhabil})
Let $K$ be an algebraically closed field of characteristic $p\geq 0$
and let $A$ be an algebra of dihedral type with two simple module.
Then $A$ is derived equivalent   to  $D(2\mathcal{B})^{k, s}(0)$ with $k\geq
s\geq 1$ or ($p=2$ and   $D(2\mathcal{B})^{k, s}(1)$ with $k\geq s\geq 1$).
\end{Prop}

\begin{Prop} \label{stableclassDiherdaltwosimples}
Let $K$ be an algebraically closed field of characteristic $p\geq
0$ and let $A$ be an algebra of dihedral type with two simple
module. Then $A$ is stably  equivalent of Morita type  to one and
exactly one of the following algebras: $D(2\mathcal{B})^{k, s}(0)$
with $k\geq s\geq 1$ or if $p=2$, $D(2\mathcal{B})^{k, s}(1)$ with
$k\geq s\geq 1$.
\end{Prop}

\begin{Proof} By the result of Pogorzar\l y (\cite[Theorem 7.3]{Pogorzaly}),
in case of characteristic two,  the algebras
$D(2\mathcal{B})^{k, s}(0)$ and $D(2\mathcal{B})^{k, s}(1)$
are not stably equivalent of Morita type.

Now for any characteristic $p$ and for  different parameters $(k,
s)\neq (k', s')$ such $k\geq s\geq 1$ and $k'\geq s'\geq 1$, if
$D(2\mathcal{B})^{k, s}(c)$ is stably equivalent to
$D(2\mathcal{B})^{k', s'}(c)$, then comparing the dimension of the
centre modulo the Reynolds ideal gives $k+s=k'+s'$. Since the
absolute values of the  determinants of the Cartan matrices are
the same, we get $ks=k's'$. This implies that $k=k'$ and $s=s'$.

\end{Proof}

\subsection{Three simple modules}

\begin{Prop} \label{stableclassDiherdalthreesimples}
Let $K$ be an algebraically closed field of characteristic
$p\geq 0$ and let $A$ be an algebra of dihedral type with two simple modules.
Then $A$ is stably  equivalent of Morita type  to one and exactly one of
the following algebras:  $D(3\mathcal{K})^{a, b, c}$ with $a\geq b\geq c\geq 1$
or   $D(3\mathcal{R})^{k, s, t, u}$ with $  s\geq t\geq u\geq k\geq 1$ and $t\geq 2$.
\end{Prop}

\begin{Proof} Holm shows \cite[page 58]{Holmhabil} that
the stable Auslander-Reiten quivers of algebras of type $D(3\mathcal{K})^{a, b, c}$
and of algebras of type $D(3\mathcal{R})^{k, s, t, u}$ is different. Hence
algebras of these two types cannot be stably equivalent of Morita type.

Again, we consider different parameters of type
$D(3\mathcal{K})^{a, b, c}$ or  of type $D(3\mathcal{R})^{k, s, t,
u}$.  By Theorem~\ref{Reynolds}, one can use  the algebra
structure of the centre modulo the Reynolds ideal to distinguish
stable equivalences classes of Morita type.

Using the the explicit basis of the centres (Holm\cite[Lemma
2.3.16]{Holmhabil}) allows to determine the the quotient
$$Z(D(3\mathcal{K})^{a, b, c})/R(D(3\mathcal{K})^{a, b, c})\simeq
K[A,B,C]/(A^a,B^b,C^c,AB,AC,BC)$$ and hence two algebras of type
$D(3\mathcal{K}^{a,b,c})$ can only be stably equivalent of Morita
type if the parameters $a,b,c$ are equal (cf
Theorem~\ref{Reynolds}).

Using (Holm\cite[Lemma 2.3.17]{Holmhabil}), we  also get that
$$Z(D(3\mathcal{R})^{k, s, t, u})/R(D(3\mathcal{R})^{k, s, t, u})\simeq
K[U,V,W,T]/(U^s,V^t,W^u,T^k,UV,UW,UT,VW,VT,WT)$$
and again two algebras of type $D(3\mathcal{R})^{k, s, t, u}$ can only be stably
equivalent of Morita type if the parameters coincide.
\end{Proof}

Although our above result is only a complete classification up to
a scalar problem in one simple module case, we can prove
nevertheless the following special case of the Auslander-Reiten
conjecture.

\begin{Cor} \label{ARconjdihedral}
Let $A$ be an indecomposable algebra which is stably
equivalent of Morita type to an algebra of dihedral type. Then
this algebra has the same number of simple modules as the algebra
of dihedral type.
\end{Cor}

\begin{Proof} By  Proposition~\ref{onlyinternal}, $A$
is necessarily of dihedral type. Then apply our classification
results above.
Notice that although we cannot determine
whether $D(1\mathcal{A})^k(0) $ and  $D(1\mathcal{A})^k(1)$ are stably equivalent
of Morita type or not, they have the same number of
simple modules.
\end{Proof}

\section{Centres of semi-dihedral and quaternion type algebras}
\label{centreSection}

We shall study the centres and the stable centres of the involved
algebras.

\subsection{Semi-dihedral type}
\label{semidihedralsection}

An algebra of semi-dihedral type with one simple module is Morita
equivalent to $SD(1\mathcal{A})_1^k$ with $k\geq 2$ or to (in
case of characteristic $2$) $SD(1\mathcal{A})_2^k(c, d)$ with
$k\geq 2$ and $(c, d)\neq (0, 0)$. Recall from \cite[Corollary
III.1.3]{ErdmannLNM} that for each of these algebras, the centre
has dimension $4k$ and  the dimension of the centre
is $k+3$. Denote by $A$ one of the above algebras. The centre
  $Z(A)$   has a $K$-basis given by
$$\{1; (XY)^i+(YX)^i; (XY)^k;  X(YX)^{k-1}; (YX)^{k-1}Y \;|\; 1\leq i\leq k-1 \}$$

\begin{Lem}\label{centrelocalsemid}
Let $K$ be an algebraically closed field and let
$A$ be one of the algebras $SD(1\mathcal{A})_1^k$ with $k\geq 2$ or (in
case of characteristic $2$) $SD(1\mathcal{A})_2^k(c, d)$ with
$k\geq 2$ and $(c, d)\neq (0, 0)$.

If $K$ is of characteristic $2$, then
$$Z(A)\simeq K[U,T,V,W]/(U^k,T^2,V^2,W^2,UT,UV,UW,TV,TW,VW)$$
and $R(A)=Z(A)\cap \mathrm{Soc}(A)=K\cdot T$.

If $K$ is of characteristic different from $2$, then
$$Z(A)\simeq K[U,V,W]/(U^{k+1},V^2,W^2,UV,UW,VW)$$
and $R(A)=Z(A)\cap \mathrm{Soc}(A)=K\cdot U^k$.
\end{Lem}

\begin{Proof} We need to identify $U$ with $XY+YX$, observe that
$((XY)+(YX))^i=(XY)^i+(YX)^i$, and identify $T$ with $(XY)^k$ and
$V$ and $W$ with the other two remaining elements. If $K$ is of
characteristic $2$, then $U^k=(XY)^k+(YX)^k=0$, and if $K$ is of
characteristic different from $2$, then
$U^k=(XY)^k+(YX)^k=2(XY)^k\neq 0$.
\end{Proof}

The Cartan matrix of the algebra $A$ is the matrix $(4k)$ of size
$1\times 1$. Recall that the dimension of the projective centre is
the $p$-rank of the Cartan matrix where $p$ is the characteristic
of the base field.   If the characteristic of $K$ divides $4k$,
then the $p$-rank of the Cartan matrix is $0$, and is $1$
otherwise.

\begin{Rem}\rm  Using the dimension of the center modulo the Reynolds ideal, we see that different values of $k$ give
different stable equivalent classes of Morita type for the above
algebras.

\end{Rem}

\medskip

Now we turn to the cases of two simple modules. An algebra of
semi-dihedral type with two  simple modules is derived equivalent
to $SD(2\mathcal{B})_1^{k,s}(c)$ with $k\geq 1$, $s\geq 2$ and
$c\in\{0,1\}$ or to $SD(2\mathcal{B})_2^{k,s}(c)$  with $k\geq 1$,
$s\geq 2$, $k+s\geq 4$ and $c\in\{0,1\}$.

\begin{Lem}\label{CenterSemidihedralTwo}
Let $A$ be the algebra  $SD(2\mathcal{B})_1^{k,s}(c)$ or the
algebra $SD(2\mathcal{}B)_2^{k,s}(c)$.

(1) If $K$ is of characteristic $2$,  then
$$Z(A) \simeq
K[u, v, w, t]/(u^k-v^s, w^2, t^2, uv,uw,vw,tw,ut,vt)$$ and
$R(A)=K\cdot u^k\oplus K\cdot w$.

(2) If $K$ is of characteristic different from $2$, then
$$Z(A)\simeq
 K[u, v, t]/(u^{k+1}, v^{s+1},  t^2, uv, ut,vt)$$
and $R(A)=K\cdot u^k\oplus K\cdot v^s$.
\end{Lem}

\begin{Proof} By \cite[IX 1.2 LEMMA]{ErdmannLNM}, a basis of the centre  of
$SD(2\mathcal{B})_1^{k,s}(c)$ is given by
$$\{1;(\alpha\beta\gamma)^i+(\beta\gamma\alpha)^i+(\gamma\alpha\beta)^i;
(\beta\gamma\alpha)^{k-1}\beta\gamma;
(\alpha\beta\gamma)^k;  \eta^j |\;1\leq i\leq k-1; 1\leq j\leq s\}$$
Now let
$$u=\alpha\beta\gamma+\beta\gamma\alpha+\gamma\alpha\beta,
v=\eta,
t=(\beta\gamma\alpha)^{k-1}\beta\gamma,
w=(\alpha\beta\gamma)^k.$$
If $\Char (K)=2$, then $u^t=v^s$; otherwise, $u^t=v^s+2w$. Hence, $w$
may be eliminated from the relations by the equation $u^k=v^s+2w$
in case  $\Char (K)\neq 2$. It is easy to verify all other relations.
An argument of comparing dimensions gives the result.

As for $SD(2B)_2^{k,s}(c)$, by \cite[IX 1.2 LEMMA]{ErdmannLNM}, a
basis of the centre of $SD(2B)_1^{k,s}(c)$ is given by
$$\{1;
(\alpha\beta\gamma)^i+(\beta\gamma\alpha)^i+(\gamma\alpha\beta)^i; (\beta\gamma\alpha)^{k-1}\beta\gamma;
(\alpha\beta\gamma)^k; \eta+(\alpha\beta\gamma)^{k-1}\alpha;
\eta^j\;|\;
1\leq i\leq k-1; 2\leq j\leq s\}$$
Now let
$$u=\alpha\beta\gamma+\beta\gamma\alpha+\gamma\alpha\beta,
v=\eta+(\alpha\beta\gamma)^{k-1}\alpha,
t=(\beta\gamma\alpha)^{k-1}\beta\gamma, w=(\alpha\beta\gamma)^k.$$
 Similar argument as above gives the result.

\end{Proof}

It is important to know that in this presentation the element $t$
is not in the socle of $SD(2\mathcal{B})_1^{k,s}(c)$ and can
therefore not be in the projective centre (cf
Proposition~\ref{Higman}).

The Cartan matrix of $SD(2\mathcal{\mathcal{B}})_1^{k,s}(c)$ with
$k\geq 1$, $s\geq 2$ and $c\in\{0, 1\}$ and of $SD(2B)_2^{k,s}(c)$
with $k\geq 1$, $s\geq 2$, $k+s\geq 4$ and $c\in\{0, 1\}$ is
$$\left(\begin{array}{cc}4k&2k\\ 2k&s+k\end{array}\right)$$
The determinant is $4ks$. If the base field is of characteristic
$2$, then the $2$-rank is $1$ if and only if $k+s$ is odd, $0$
else. If the base field is of characteristic $p>2$, then the
$p$-rank of the Cartan matrix is $1$ if and only if $p$ divides
$k$ or $p$ divides $s$ but not both;  the $p$-rank is $0$ if and
only if $p$ divides $k$ and $s$; and the $p$-rank is $2$ if and
only if $p$ divides neither $k$ nor $s$. If the characteristic of
the field is $0$, then the rank of the Cartan matrix is $2$.

\bigskip

Recall that Holm proved in \cite[Lemma 2.4.16]{Holmhabil} that the
centre of $SD(3\mathcal{K})^{a,b,c}$  with $a\geq b\geq c\geq 1$
and $a\geq 2$ has a basis given by
$$\{1, (\beta\gamma+\gamma\beta)^{i_1};
(\kappa\lambda+\lambda\kappa)^{i_2};
(\delta\eta+\eta\delta)^{i_3}; (\lambda\kappa)^b; (\beta\gamma)^b;
(\delta\eta)^c\;|\;
1\leq i_1< a, 1\leq i_2<b ,1\leq i_3\leq c \}$$
and so

\begin{Lem}\label{CenterSemidihedralThree}
\begin{eqnarray*}
Z(SD(3\mathcal{K})^{a,b,c})&\simeq&
K[A,B,C,S_1,S_2,S_3]/(A^{a+1},B^{b+1},C^{c+1},A^a-S_2-S_3,B^b-S_3-S_1,\\
&&C^c-S_1-S_2,AS_i,BS_i,CS_i,S_iS_j,AB,AC,BC; i,j\in\{1,2,3\})
\end{eqnarray*}
and $R(SD(3\mathcal{K})^{a,b,c})=K\cdot A^a\oplus  K\cdot
B^b\oplus K\cdot C^c$.

\end{Lem}

\begin{Proof} Let $$A=\beta\gamma+\gamma\beta,
B=\kappa\lambda+\lambda\kappa, C=\delta\eta+\eta\delta,
S_1=\lambda\beta\delta, S_2=\delta\lambda\beta,
S_3=\beta\delta\lambda.$$
Then it is a straight forward verification that
$A, B, C, S_1, S_2, S_3$ satisfy the relations on the right-handed
side. Now the isomorphism follows from a dimension argument.
\end{Proof}

\begin{Cor} Let $a\geq b\geq c\geq 1$ with $a\geq 2$ and let
$a'\geq b'\geq c'\geq 1$ with $a'\geq
2$. Then $SD(3\mathcal{K})^{a,b,c}$ is stably equivalent of Morita
type to $SD(3\mathcal{K})^{a',b',c'}$ if and only if $a=a'$,
$b=b'$ and $c=c'$.
\end{Cor}

\begin{Proof} As in the proof of
Proposition~\ref{stableclassDiherdalthreesimples}, one can
consider the centre modulo the Reynolds ideal. The Reynolds ideal
of the centre is of dimension three and is spanned by the elements
$S_1,S_2,S_3$. Hence
$$Z(SD(3\mathcal{K})^{a,b,c})/R(SD(3\mathcal{K})^{a,b,c})\simeq
K[A,B,C]/(A^a,B^b,C^c,AB,AC,BC)$$
so that an isomorphism of the centres modulo the Reynolds ideal
implies that the parameters are identical. The statement then follows from
Theorem~\ref{Reynolds}.

\end{Proof}

The Cartan matrix of $SD(3\mathcal{K})^{a,b,c}$ equals (cf
\cite{ErdmannLNM})
$$\left(\begin{array}{ccc}a+b&a&b\\ a&a+c&c\\ b&c&b+c\end{array}\right)$$
which has determinant $4abc$.
Since all
coefficients of the Cartan matrix  are positive integers,  the
rank of the Cartan matrix for fields of characteristic $0$ is
always $3$.

Suppose that $K$ is a base field of characteristic
$p>2$. \\
The $p$-rank $0$ occurs if and only if all parameters
$a,b,c$ are divisible by $p$; \\
the $p$-rank $1$ occurs if and only
if exactly one parameter $a,b,c$ is not divisible by $p$; \\
the
$p$-rank is $2$ if and only if exactly one of the parameters is
divisible by $p$; \\
the $p$-rank is $3$ if and only if $p$ doesn't
divide $abc$.

Suppose that $K$ is a base field of characteristic
$2$. \\
The $2$-rank $0$ occurs if and only if all parameters $a,b,c$
are all even;\\
the $2$-rank $1$ occurs if and only if exactly one
parameter $a,b,c$ is odd; \\
the $2$-rank is $2$ if and only if at
least two of $a, b, c$ are odd.

\subsection{Quaternion type}
\label{quaternionsection}

An algebra of  quaternion type with one simple module is Morita
equivalent to $Q(1\mathcal{A})_1^{\ell}$ with ${\ell}\geq 2$ or to
$Q(1\mathcal{B})_2^{\ell}(c,d)$ with ${\ell}\geq 2$ and $(c,
d)\neq (0, 0)$.
 Again, by \cite{ErdmannLNM} the centre is of dimension $\ell+3$
and the algebra is of dimension $4\ell$. Let $A$ be one of the
above algebras. In the above presentation, the centre has a
$K$-basis given by
$$\{1; (XY)^i+(YX)^i; (XY)^\ell;
X(YX)^{\ell-1}; (YX)^{\ell-1}Y \;|\;  1\leq i\leq \ell-1\}.$$

\begin{Lem}
(1) If $K$ is of characteristic $2$, then
$$Z(A)\simeq K[U,T,V,W]/(U^k,T^2,V^2,W^2,UT,UV,UW,TV,TW,VW)$$
and $R(A)=Z(A)\cap soc(A)=K\cdot T$.

(2) If $K$ is of characteristic different from $2$, then
$$Z(A)\simeq K[U,V,W]/(U^{k+1},V^2,W^2,UV,UW,VW)$$
and $R(A)=Z(A)\cap soc(A)=K\cdot U^k$.
\end{Lem}

\begin{Proof}
The proof is a straight forward verification.
\end{Proof}

\bigskip

An algebra of quaternion type with two simple modules is derived
equivalent to $Q(2\mathcal{B})_1^{k,s}(a, c)$ with $k\geq 1$,
$s\geq 3$ and $a\neq 0$. By \cite[IX 1.2 LEMMA]{ErdmannLNM}, the
centre of $Q(2\mathcal{B})_1^{k,s}(a, c)$ has a basis
$$
\{1;
(\alpha\beta\gamma)^{i}+(\beta\gamma\alpha)^{i}+(\gamma\alpha\beta)^{i},
(\beta\gamma\alpha)^{k-1}\beta\gamma,
(\alpha\beta\gamma)^k,
\eta+(\alpha\beta\gamma)^{k-1}\alpha,
\eta^j\;|\;  1\le i\le k-1;2\le j\le s\}.
$$
By a similar proof as that of
Proposition~\ref{CenterSemidihedralTwo}, we have

\begin{Lem}\label{CenterQuatTwo}
(1) If $\Char(K) =2$, then $$Z(Q(2\mathcal{B})_1^{k,s}(a, c))
\simeq K[u, v, w, t]/(u^k-v^s, w^2, t^2, uv,uw,vw,tw,ut,vt)$$ and
$R(A)=Z(A)\cap soc(A)=K\cdot u^k\oplus K\cdot w$.

(2) If $\Char(K) \neq 2$, then $$Z(Q(2\mathcal{B})_1^{k,s}(a, c))
\simeq K[u, v, t]/(u^{k+1}, v^{s+1},    t^2, uv, ut,vt)$$ and
$R(A)=Z(A)\cap soc(A)=K\cdot u^k\oplus K\cdot v^s$.
\end{Lem}

The Cartan matrix of $Q(2\mathcal{B})_1^{k,s}(a, c)$  is
$$\left(\begin{array}{cc}4k&2k\\ 2k&k+s\end{array}\right).$$

\bigskip

 An algebra of quaternion type with three simple modules is
derived equivalent to $Q(3\mathcal{K})^{a,b,c}$ with $a\geq b\geq
c\geq 1$, $b\geq 2$ and $(a, b, c)\neq (2, 2, 1)$ or to
$Q(3\mathcal{A})_1^{2,2}(d)$ with $d\not\in\{ 0, 1\}$.

The dimension of the centre of $Q(3\mathcal{K})^{a,b,c}$ is $a+b+c+1$ and
has a basis
$$\{1, (\beta\gamma+\gamma\beta)^{i_1},
(\kappa\lambda+\lambda\kappa)^{i_2},
(\delta\eta+\eta\delta)^{i_3}; (\lambda\kappa)^b; (\beta\gamma)^b;
(\delta\eta)^c\;|\; 1\leq i_1< a; 1\leq i_2<b; 1\leq i_3\leq c \} .$$
The Cartan matrix of the algebra $Q(3\mathcal{K})^{a,b,c}$ is
$$\left(\begin{array}{ccc}a+b&a&b\\ a&a+c&c\\ b&c&b+c\end{array}\right)$$

The dimension of the centre of  $Q(3\mathcal{A})_1^{2,2}(d)$ is
$6$ and has a basis
$$\{1, \beta\gamma+\gamma\beta+d\eta\delta,
\beta\gamma+\eta\delta+\delta\eta, (\beta\gamma)^2,
(\gamma\beta)^2=d(\delta\eta)^2,  (\eta\delta)^2  \}.$$
The Cartan matrix of  $Q(3\mathcal{A})_1^{2,2}(d)$ is
$$\left(\begin{array}{ccc}4&2&2\\ 2&3&1\\ 2&1&3\end{array}\right)$$
Indeed, the fact that the above elements are central as is readily
verified and the dimensions are as they should be. The statement
on the Cartan matrix is taken from \cite{ErdmannLNM}.

\begin{Lem}
We have
\begin{eqnarray*}
Z(Q(3\mathcal{K})^{a,b,c})&\simeq&
K[A,B,C,S_1,S_2,S_3]/(A^{a+1},B^{b+1},C^{c+1},A^a-S_2-S_3,B^b-S_3-S_1,\\
&&C^c-S_1-S_2,AS_i,BS_i,CS_i,S_iS_j,AB,AC,BC; i,j\in\{1,2,3\})\\
Z(Q(3\mathcal{A})_1^{2,2}(d))&\simeq&
K[A,B,C,S_1,S_2,S_3]/(A^{3},B^{3},C^{2},A^2-S_2-S_3,B^2-S_3-S_1,\\
&&C-S_1-S_2,AS_i,BS_i,CS_i,S_iS_j,AB,AC,BC;
i,j\in\{1,2,3\})\end{eqnarray*}
\end{Lem}

\begin{Proof} The proof for $Q(3\mathcal{K})^{a,b,c}$ is identical to the
one of Lemma~\ref{CenterSemidihedralThree}. For
$Q(3\mathcal{A})_1^{2,2}(d)$, let
\begin{eqnarray*}A&:=&\beta\gamma+\gamma\beta+d\eta\delta,\\
B&:=&\beta\gamma+\eta\delta+\delta\eta,\\
C&:=&(1-d)(\delta\eta)^2+d^2(\eta\delta)^2,\\
S_1&:=&(1-d)(\delta\eta)^2,\\
S_2&:=&d^2(\eta\delta)^2,\\
S_3&:=&(\beta\gamma)^2+d(\delta\eta)^2.
\end{eqnarray*}
The rest is a straight forward verification.
\end{Proof}

Note that, in order to simplify the notation we may put
$Q(3\mathcal{K})^{2,2,1}=Q(3\mathcal{A})_1^{2,2}(d)$.

\begin{Cor} Let $a\geq b\geq c\geq 1$ with $b\geq 2$ and let
$a'\geq b'\geq c'\geq 1$ with $b'\geq
2$.
Then $Q(3\mathcal{K})^{a,b,c}$ is stably equivalent of Morita
type to $Q(3\mathcal{K})^{a',b',c'}$ if and only if $a=a'$,
$b=b'$ and $c=c'$.
\end{Cor}

\begin{Proof} As in the proof of
Proposition~\ref{stableclassDiherdalthreesimples}, one can
consider the centre modulo the Reynolds ideal. The Reynolds ideal
of the centre is of dimension three and is spanned by the elements
$S_1,S_2,S_3$ and as in the proof of the semi-dihedral type
$$Z(Q(3\mathcal{K})^{a,b,c})/R(Q(3\mathcal{K})^{a,b,c})\simeq
K[A,B,C]/(A^a,B^b,C^c,AB,AC,BC)$$
so that as in the semi-dihedral case an isomorphism
of the centres modulo the Reynolds ideal
implies that the parameters are identical. Finally apply
Theorem~\ref{Reynolds}.
\end{Proof}

\section{Algebras with stable centres and Cartan data as for
semi-dihedral and quaternion type; stable equivalences}
\label{thecentralpart}

{\bf For $a, b, c\geq 1$, let $A_3^{a,b,c}$} be a basic indecomposable
symmetric $K$-algebra with centre isomorphic to
\begin{eqnarray*}
Z(A_3^{a,b,c})&\simeq&
K[A,B,C,S_1,S_2,S_3]/(A^{a+1},B^{b+1},C^{c+1},A^a-S_2-S_3,B^b-S_3-S_1,\\
&&C^c-S_1-S_2,AS_i,BS_i,CS_i,S_iS_j,AB,AC,BC; i,j\in\{1,2,3\})
\end{eqnarray*}
and Cartan matrix
$$\left(\begin{array}{ccc}a+b&a&b\\ a&a+c&c\\ b&c&b+c\end{array}\right)$$
and the Reynolds ideal $R(A_3^{a,b,c})=KA^a\oplus KB^b\oplus
KC^c$.

{\bf For $k, s\geq 1$, let $A_2^{k,s}$} be a basic indecomposable
symmetric algebra  with Cartan matrix
$$\left(\begin{array}{cc}4k&2k\\ 2k&s+k\end{array}\right)$$
so that in case $K$ is of characteristic $2$,
the centre
$$Z(A_2^{k,s})\simeq
 K[u, v, w, t]/(u^k-v^s, w^2, t^2, uv,uw,vw,tw,ut,vt)$$
and   $R(A_2^{k,s})= K u^k\oplus K w$ and if if $K$ is of
characteristic different from $2$, then $$Z(A_2^{k,s})\simeq
 K[u, v, t]/(u^{k+1}, v^{s+1},   t^2, uv, ut,vt)$$
and   $R(A_2^{k,s})= K u^k\oplus K v^s$

{\bf For $\ell\geq 2$, let $A_1^\ell$} be a basic indecomposable
symmetric algebra of dimension $4\ell$ so that in case $K$ is of
characteristic $2$,
$$Z(A_1^\ell)\simeq K[U,T,V,W]/(U^\ell,T^2,V^2,W^2,UT,UV,UW,TV,TW,VW)$$
and the Reynolds ideal $R(A_1^\ell)=K\cdot T$ and if $K$ is of
characteristic different from $2$, then
$$Z(A_1^\ell)\simeq K[U,V,W]/(U^{\ell+1},V^2,W^2,UV,UW,VW)$$
and $R(A_1^\ell)=K\cdot U^\ell$.

\subsection{Two simples versus three simples; characteristic different from $2$}

Concerning the relations of $Z(A_3^{a,b,c})$ we see that the
elements $S_1+S_2, S_2+S_3, S_3+S_1$ generate the same space as
$S_1,S_2,S_3$, which is the whole socle of the algebra, if and
only if
$$\left(\begin{array}{ccc}1&1&0\\0&1&1\\1&0&1\end{array}\right)$$
is a regular matrix. This is the case if and only if $K$ is a
field of characteristic different from $2$. Therefore if the
characteristic of the base field is different from $2$, we get
$$Z(A_3^{a,b,c})\simeq K[A,B,C]/(A^{a+1},B^{b+1}, C^{c+1},AB,AC,BC).$$

\begin{Lem}\label{2vers3carpiii}
Let $K$ be an algebraically closed field of characteristic $p>2$
or of characteristic $0$. Suppose $p\not| (ks)$. Then $A_2^{k,s}$
and $A_3^{a,b,c}$ cannot be stably equivalent of Morita type.
\end{Lem}

\begin{Proof}
If $K$ is of characteristic $0$, then the rank of the Cartan
matrix of $A_3^{a,b,c}$ is $3$ and the rank of the Cartan matrix
of $A_2^{k,s}$ is $2$. Hence
\begin{eqnarray*}
Z^{st}(A_2^{k,s})&=&K[u,v,t]/(u^k,v^s,t^2,uv,ut,vt)\\
Z^{st}(A_3^{a,b,c})&=&K[A,B,C]/(A^a,B^b,C^c,AB,AC,BC)
\end{eqnarray*}
The same holds if $p\not| (ks)$ because then $abc=ks$ since the
Cartan determinants coincide, and since therefore the Cartan
matrices are regular. Hence in order to get the stable centres
isomorphic we may assume $c=2$, $a=k$, $b=s$, else a permutation
of the letters $a,b,c$ will do. Hence $ks=abc$ becomes $ks=2ks$, a
contradiction.
\end{Proof}

\begin{Lem}\label{2vers3carpi}
Let $K$ be an algebraically closed field of characteristic $p>2$.
Suppose $p|k$ and $p|s$. Then $A_2^{k,s}$ and $A_3^{a,b,c}$ cannot
be stably equivalent of Morita type.
\end{Lem}

\begin{Proof}
Suppose that the algebras are stably equivalent of Morita type. We
know that the Cartan determinants coincide and hence $p^2| abc$.
If $p$ divides $a$ and $b$ and $c$, then the Cartan matrix of
$A_3^{a,b,c}$ three elementary divisors divisible by
$p$, which implies that the stable Grothendieck group of
$A_3^{a,b,c}$ tensored by $K$ has rank $3$.
This gives a contradiction since the
stable Grothendieck group of $A_2^{k,s}$ tensored by $K$ can only be of rank $2$
at most.

Since the $p$-rank of the Cartan matrix of $A_2^{k,s}$
is $0$, we get that
\begin{eqnarray*}
Z^{st}(A_2^{k,s})&=&K[u,v,t]/(u^{k+1},v^{s+1},t^2,uv,ut,vt)
\end{eqnarray*}

If $p$ divides two of the parameters $a$, $b$, $c$, then
\begin{eqnarray*}
Z^{st}(A_3^{a,b,c})&=&K[A,B,C]/(A^{a+1},B^{b+1},C^{c+1},AB,AC,BC,
\lambda_AA^a+\lambda_BB^b+\lambda_CC^c)
\end{eqnarray*}
for some parameters $\lambda_A,\lambda_B,\lambda_C$ not all $0$.
If $p$ divides only one of the parameters $a$, $b$, $c$, then
\begin{eqnarray*}
Z^{st}(A_3^{a,b,c})&=K[A,B,C]/&(A^{a+1},B^{b+1},C^{c+1},AB,AC,BC,\\
&&\lambda_AA^a+\lambda_BB^b+\lambda_CC^c,\mu_AA^a+\mu_BB^b+\mu_CC^c)
\end{eqnarray*}
for a matrix
$$\left(\begin{array}{ccc}\lambda_A&\lambda_B&\lambda_C\\
\mu_A&\mu_B&\mu_C\end{array}\right)$$ of rank $2$.

The socle of the stable centre of $A_2^{k,s}$ is
three-dimensional, and so we need to assure that this is the case
of the stable centre of $A_3^{a,b,c}$ as well. But this implies
that the projective centre of $A_3^{a,b,c}$ is generated by $B^b$
and $C^c$ (say) if $p$ divides only one of $a,b,c$ and by $C^c$
(say) if $p$ divides two of the parameters $a$, $b$ and $c$.

In the first case, $p$ divides only one of the parameters $a,b,c$,
we get $\{a+1,b,c\}=\{k+1,s+1,2\}$, taken with multiplicities. If
$c=2$ (or $b=2$, case which is studied analogously), then $abc=ks$
becomes $2ab=ks$. Moreover, $a\in\{k,s\}$ and $b\in\{k+1,s+1\}$
gives $2k(s+1)=ks$ or $2s(k+1)=ks$. Hence $ks+k=0$ or $ks+s=0$,  a
contradiction. Hence $a=1$. But then $b=k+1$ and $c=s+1$. Now, $p$ was
assumed to divide $k$ and $s$, and so $p$ divides none of the
parameters $a,b,c$, a contradiction to the hypothesis.

If $p$ divides two of the parameters $a$, $b$, $c$, the projective
centre is one-dimensional and we got that $C^c$, say, generates
the projective centre. Then $\{a+1,b+1,c\}=\{k+1,s+1,2\}$ again
taken with multiplicities. If
$c=2$, then $k=a$ and $b=s$, say, and $abc=2ks=ks$, a
contradiction. Hence by symmetry we may assume $b=1$. If
$c=s+1$ and $k=a$, then
$abc=k(s+1)=ks$ gives a contradiction; if $c=k+1$ and $a=s$ the
same contradiction holds.
\end{Proof}

\begin{Lem}\label{2vers3carpii}
Let $K$ be an algebraically closed field of characteristic $p>2$.
Suppose $p|k$ and $p\not|s$ or $p|s$ and $p\not|k$. Then
$A_2^{k,s}$ and $A_3^{a,b,c}$ cannot be stably equivalent of
Morita type.
\end{Lem}

\begin{Proof}
The hypothesis implies that the projective centre of $A_2^{k,s}$
is one-dimensional, and hence there are parameters $\nu_u,\nu_v$
not both $0$ with
\begin{eqnarray*}
Z^{st}(A_2^{k,s})&=&K[u,v,t]/(u^{k+1},v^{s+1},t^2,\nu_uu^k+\nu_vv^s,uv,ut,vt)
\end{eqnarray*}
Again, as before, the stable Grothendieck groups need to be
isomorphic and so not all parameters $a$, $b$, $c$ can be
divisible by $p$. Actually, since one of the parameters $k$ and
$s$ is not divisible by $p$, one of the elementary divisors of the
Cartan matrix of $A_2^{k,s}$ is not divisible by $p$ and the other
is divisible by $p$. Hence one of the elementary divisors of the
Cartan matrix of $A_3^{a,b,c}$ is $1$, one is not divisible by $p$
and the third is divisible by $p$. If $p$ divides two of the
parameters $a$, $b$ and $c$, then two elementary divisors of the
Cartan matrix of $A_3^{a,b,c}$ are divisible by $p$, whence a
contradiction. Hence $p$ divides exactly one of the parameters
$a,b,c$, and the projective centre of $A_3^{a,b,c}$ is
two-dimensional. We get
\begin{eqnarray*}
Z^{st}(A_3^{a,b,c})&=K[A,B,C]/&(A^{a+1},B^{b+1},C^{c+1},AB,AC,BC,\\
&&\lambda_AA^a+\lambda_BB^b+\lambda_CC^c,\mu_AA^a+\mu_BB^b+\mu_CC^c)
\end{eqnarray*}
for a matrix
$$\left(\begin{array}{ccc}\lambda_A&\lambda_B&\lambda_C\\
\mu_A&\mu_B&\mu_C\end{array}\right)$$ of rank $2$.

Suppose $\nu_u=0$ or $\nu_v=0$. By symmetry we may suppose
$\nu_u=0$. Then the socle of $Z^{st}(A_2^{k,s})$ is
three-dimensional. Hence in order to get this we need to have that
$B^b$ and $C^c$, say, generate the projective centre of
$A_3^{a,b,c}$. But then $\{a+1,b,c\}=\{k,s+1,2\}$, taken with
multiplicities. The case $a=1$ gives a contradiction to $abc=ks$
as well as the case $c=2$ (or likewise $b=2$).

Hence $\nu_u\neq 0\neq \nu_v$. The socle of $A_2^{k,s}$ is
two-dimensional and therefore the socle of
$Z^{st}(A_3^{a,b,c})$ has to be two-dimensional as well. This
implies that one of the elements $A^a, B^b, C^c$ has to be in the
projective centre, say $C^c$. Therefore
\begin{eqnarray*}
Z^{st}(A_3^{a,b,c})&=&K[A,B,C]/(A^{a+1},B^{b+1},C^{c},AB,AC,BC,
\lambda_AA^a+\lambda_BB^b)
\end{eqnarray*}
This give $c=2$ and $\{a,b\}=\{k,s\}$. Now, the equality of Cartan
determinants $abc=ks$ is not satisfied, a contradiction.
\end{Proof}

\subsection{Two simples versus three simples; characteristic $2$}

We are now dealing with the case $p=2$. Recall that
$$Z(A_2^{k,s})\simeq K[u, v, w, t]/(u^k-v^s, w^2, t^2, uw,vw,tw,ut,vt)$$
and
\begin{eqnarray*}
Z(A_3^{a,b,c})&\simeq& K[A,B,C,S_1,S_2,S_3]/
(A^{a+1},B^{b+1},C^{c+1},A^a-S_1-S_2,B^b-S_2-S_3,\\
&&C^c-S_1-S_3,AS_i,BS_i,CS_i,S_iS_j,AB,AC,BC; i,j\in\{1,2,3\})
\end{eqnarray*}
In case $p=2$ the subspace of the socle of the algebra generated
by $S_1+S_2,S_2+S_3,S_1+S_3$ is of codimension $1$, namely given
by the condition
$$(S_1+S_2)+(S_2+S_3)+(S_1+S_3)=0$$
and so
$$A^a+B^b+C^c=0.$$
Hence, in characteristic $2$ we get
\begin{eqnarray*}
Z(A_3^{a,b,c})&\simeq&
K[A,B,C,S]/(A^{a+1},B^{b+1},C^{c+1},S^2,A^a+B^b+C^c,AS,BS,CS,AB,AC,BC)
\end{eqnarray*}

We shall show

\begin{Prop}\label{2vers3car2}
Let $K$ be an algebraically closed field of characteristic $2$.
$A_3^{a,b,c}$ cannot be   stably equivalent of Morita type to
$A_2^{k,s}$.
\end{Prop}

The proof will consist of two technical lemmata~\ref{paritiesofparameters}
and \ref{kplussOdd} contradicting each other.

\begin{Lem} \label{paritiesofparameters}
Let $K$ be an algebraically closed field of characteristic $2$.
Suppose that  $A_3^{a,b,c}$ is stably equivalent of Morita type to
$A_2^{k,s}$. Then $k$ and $s$ are both even, two of $a, b, c$ are
odd and the third is even.
\end{Lem}

\begin{Proof}
Supposing that  $A_3^{a,b,c}$ is stably equivalent of Morita type
to $A_2^{k,s}$, then $abc=ks$.
The dimensions of the centre modulo the Reynolds ideal gives
$k+s=a+b+c-2$.

If $a$, $b$ and $c$ are all even, then all elementary divisors of
the Cartan matrix of the algebra $A_3^{abc}$ are divisible by $2$
and hence the stable Grothendieck group tensored by $K$ of
$A_3^{abc}$ is of dimension $3$. But the stable Grothendieck group
(tensored by $K$) of $A_2^{k,s}$ is $2$ at most. Therefore this
cannot happen. We get that at least one of $a$, $b$ or $c$ is odd.

If $k$ and $s$ are both odd, then $a, b, c$ are all odd, but then
$k+s=a+b+c-2$ cannot hold.


Suppose now that $ks=abc$ is even.

The stable centre of $A_2^{k,s}$ is of dimension $k+s+2$ if $k+s$
is even and of dimension $k+s+1$ if $k+s$ is odd. Hence the stable
centre of $A_2^{k,s}$ is of even dimension in any case. If two
parameters, say $a$ and $b$ are even, then the $2$-rank of the
Cartan matrix of $A_3^{a,b,c}$ is one, and hence the stable centre
of $A_3^{a,b,c}$ is of dimension $a+b+c$. Since $c$ is odd $a+b+c$
is odd and we get a contradiction.

We have proved that two among  $a$, $b$ or $c$ are odd and the third is odd
now and the equality of the dimensions of the quotient of the centres by
the Reynolds ideals $k+s=a+b+c-2$ shows that $k$ and $s$ are
both even.



\end{Proof}






\begin{Lem}\label{kplussOdd}
Let $K$ be an algebraically closed field of characteristic $2$.
Suppose that  $A_3^{a,b,c}$ is stably equivalent of Morita type to
$A_2^{k,s}$. Then $k+s$ is odd.
\end{Lem}

\begin{Proof}
If $k+s$ is even then
$$Z(A_2^{k,s})=Z^{st}(A_2^{k,s})\simeq
 K[u, v, w, t]/(u^k-v^s, w^2, t^2, uv,uw,vw,tw,ut,vt)$$
is isomorphic to the quotient of
\begin{eqnarray*}
Z(A_3^{a,b,c})&\simeq&
K[A,B,C,S]/(A^{a+1},B^{b+1},C^{c+1},S^2,A^a+B^b+C^c,AS,BS,CS,AB,AC,BC)
\end{eqnarray*}
by a two-dimensional subspace of $<A^a,B^b,S>_K$ since by
Lemma~\ref{paritiesofparameters} at most one of the parameters
$a,b,c$ is even. The socle of $A_3^{a,b,c}$ is generated by
$A^a,B^b,S$ and the projective centre is generated by two elements
$\mu_AA^a+\mu_BB^b+\mu_SS$ and $\nu_AA^a+\nu_BB^b+\nu_SS$.

{\bf If $\mu_S\neq 0$ or $\nu_S\neq 0$} then the stable centre of
$SD(3K)$ is isomorphic to
\begin{eqnarray*}
Z^{st}(A_3^{a,b,c})&\simeq&
K[A,B,C]/(A^{a+1},B^{b+1},C^{c+1},A^a+B^b+C^c,\mu_A'A^a+\mu_B'B^b,AB,AC,BC)
\end{eqnarray*}
for some parameters $\mu'_A$ and $\mu'_B$ not both $0$. In case
$min(k,s)\geq 2$ we consider the quotient modulo the radical
squared of both stable centres. The one of $A_2^{k,s}$ has a basis
given by $1,u,v,w,t$ and the one of $A_3^{a,b,c}$ is $K$-linearly generated
by $1,A,B,C$ (knowing that in case some of the
parameters $a,b,c$ are $1$, then these $4$ elements are not
linearly independent). In any case this gives a contradiction and
so $s=1$ or $k=1$ in this case. So $k$ and $s$ are both odd, which
is impossible by Lemma~\ref{paritiesofparameters}.


{\bf So we get $\mu_S=0=\nu_S$. } But then the projective centre
of $A_3^{a,b,c}$ is generated by $A^a$ and $B^b$ and we get
\begin{eqnarray*}
Z^{st}(A_3^{a,b,c})&\simeq&
K[A,B,C,S]/(A^{a},B^{b},C^{c},S^2,AS,BS,CS,AB,AC,BC)
\end{eqnarray*}
which needs to be isomorphic to
$$Z^{st}(A_2^{k,s})=K[u, v, w, t]/(u^k-v^s, w^2, t^2, uv,uw,vw,tw,ut,vt)$$
By symmetry we may assume again $a\geq b\geq c$ and $k\geq s$. If
$c\geq 2$, the socle of $Z^{st}(A_3^{a,b,c})$ is four-dimensional,
whereas the socle of $Z^{st}(A_2^{k,s})$ is three-dimensional.
Hence $c=1$. But then, comparing the quotient modulo the radical
squared gives $s=1$ and therefore $a=k+1$ and $b=2$. The equation
$abs=ks$ becomes $2(k+1)=k$, a contradiction.
\end{Proof}

\subsection{One simple versus two simples}

We shall deal with the possibility that
an algebra of type $A_1^{\ell}$ with one simple
module is stably equivalent of Morita type to an algebra of type
$A_2^{k,s}$ with $2$ simple modules.

\begin{Lem}\label{pythagoras}\label{1vers2car2}
An algebra of type $A_1^{\ell}$ with one simple module cannot be
stably equivalent of Morita type to an algebra of type
$A_2^{k,s}$.
\end{Lem}

\begin{Proof} Suppose that $A_1^\ell$ and $A_2^{k,s}$ are stably equivalent
of Morita type.
Since the Cartan determinants are equal,  we have $4\ell=4ks$.
Since the centre modulo the Reynolds ideal is invariant under a
stable equivalence of Morita type in our case by Theorem~\ref{Reynolds}, we have
$(\ell+3)-1=(k+s+2)-2$.
This means that $k$ and $s$ are integer solutions of the
equation of second order
$X^2-(\ell+2)X+\ell=0$. But
the discriminant of this equation
is equal to $(\ell+2)^2-4\ell=\ell^2+4$ which should be a square
of an integer $m$.
We look for pythagorean triples $(2,\ell,m)$. It is well known that
$\ell$ and $m$ have to be odd, whence writing down the pythagorean
equation one gets
$m=1$ or $l=1$. This contradiction proves that such a triple does not exist.
\end{Proof}

\subsection{One simple versus three simples}

\begin{Lem}\label{1vers3car2}
An algebra of type $A_1^{\ell}$  is not stably equivalent of
Morita type to an algebra of type $A_3^{a,b, c}$.
\end{Lem}

\begin{Proof}
The proof follows the lines of the proof of Lemma~\ref{pythagoras}.
The equality of Cartan determinants give
$4\ell=4abc$ and the centre modulo the Reynolds ideal give
$(\ell+3)-1=(a+b+c+1)-3$. This means that we have two equalities.
$\ell=abc $ and $\ell=a+b+c-4$. By symmetry we
may suppose $a\geq b\geq c\geq 1$.
Then $c=1$, otherwise
$$\ell=abc\geq 4a>3a\geq a+b+c-4.$$
But now
$l=ab$ and $l=a+b-3$. The same argument gives $b=1$ and now we
have $a=\ell=a-2$, which is a contradiction.
\end{Proof}

We resume the situation.

\begin{Prop}\label{result}
Let $K$ be a field of characteristic $2$.
\begin{itemize}
\item  $A_3^{a,b,c}$ cannot be  stable equivalent of Morita type
to $A_2^{k,s}$ (cf Proposition~\ref{2vers3car2}).

\item $A_2^{k, s}$ cannot be  stable equivalent of Morita type to
$A_1^{\ell}$ (cf Lemma~\ref{1vers2car2}).

 \item $A_3^{a,b,c}$ cannot be stably equivalent of
Morita type to $A_1^\ell$ (cf Lemma~\ref{1vers3car2}).
\end{itemize}
Let $K$ be a field of characteristic different from $2$.
\begin{itemize}
\item There is no stable equivalence of Morita type between
$A_3^{a,b,c}$ and $A_2^{k,s}$ (cf Lemma~\ref{2vers3carpiii},
Lemma~\ref{2vers3carpi} and Lemma~\ref{2vers3carpii}).
\item There is no stable
equivalence of Morita type between $A_1^\ell$ and $A_2^{k,s}$
(cf Lemma~\ref{1vers2car2}).
\item There is no stable equivalence of Morita type between
$A_1^\ell$ and $A_3^{a,b,c}$ (cf Lemma~\ref{1vers3car2}).
\end{itemize}
\end{Prop}

Although we cannot classify completely algebras of semi-dihedral
and quaternion type up to stable equivalences of Morita type, we
can nevertheless prove the following

\begin{Cor}\label{ARconjfortametype}
 Let $A$ be an indecomposable algebra which is stably equivalent
 to an algebra $B$ of semi-dihedral type (resp. quaternion type). Then
 $A$ has the same number of simple modules as $B$.
 \end{Cor}

 \begin{Proof}
 This is an immediate consequence of the above
 proposition.
 \end{Proof}

\section{The main theorem and concluding remarks}
\label{thestatement}

We resume the results of this paper in a single theorem. We use the notations
introduced above, which coincides with the notations in \cite{ErdmannLNM}
or \cite{Holmhabil}.

\begin{Thm}\label{main}
Let $K$ be an algebraically closed field.

Suppose $A$ and $B$ are indecomposable algebras which
are stably equivalent of Morita type.

\begin{itemize}
\item
If $A$ is an algebra of dihedral type, then $B$ is of
dihedral type. If $A$ is of semi-dihedral type, then
$B$ is of semi-dihedral type. If $A$ is of quaternion type then
$B$ is of quaternion type.

\item If $A$ and $B$ are of dihedral, semidihedral or quaternion type, then
$A$ and $B$ have the same number of simple modules.

\item Let $A$ be an algebra of
dihedral type.

\begin{enumerate}
\item
If $A$ is local, then $A$ is stably
equivalent of Morita type to one and exactly one  algebra in the
following list:
\begin{itemize}
\item $A_1(n, m)$ with $m\geq n\geq 2$ and $m+n>4$;
\item $C_1 $;
\item $D(1\mathcal{A})_1^k $ with $k\geq 2$;
\item  if $p=2$, $B_1$  and $D(1\mathcal{A})_2^k(d)$
with $k\geq 2$ and $d\in \{0, 1\}$, except that we don't
  know whether $D(1\mathcal{A})^k(0) $ and  $D(1\mathcal{A})^k(1)$ are
  stably equivalent of Morita type or not.
\end{itemize}

\item
If $A$ has two simple modules,
then $A$ is stably  equivalent of Morita type  to one and exactly
one of the following algebras:
$D(2\mathcal{B})^{k, s}(0)$ with $k\geq s\geq 1$ or if $p=2$, $D(2\mathcal{B})^{k, s}(1)$
with $k\geq s\geq 1$. 

\item
If $A$ has three simple modules then
$A$ is stably  equivalent of Morita type  to one and exactly one of
the following algebras:  $D(3\mathcal{K})^{a, b, c}$ with $a\geq b\geq c\geq 1$
or   $D(3\mathcal{R})^{k, s, t, u}$ with $  s\geq t\geq u\geq k\geq 1$ and $t\geq 2$.
\end{enumerate}

\item Let $A$ be an algebra of semi-dihedral type.
\begin{enumerate}
\item If $A$ has one simple module then $A$ is stably equivalent
of Morita type to one of the following algebras:
$SD(1\mathcal{A})_1^k$ for $k\geq 2$ or
$SD(1\mathcal{A})_2^k(c,d)$ for $k\geq 2$ and $(c,d)\neq (0,0)$ if
the characteristic of $K$ is $2$. Different parameters $k$ yield
algebras in different
 equivalence classes of Morita type.

\item
If $A$ has two simple modules then $A$ is stably equivalent of Morita type to
$SD(2\mathcal B)_1^{k,s}(c)$ for $k\geq 1,s\geq 2,c\in\{0,1\}$
or to
$SD(2\mathcal B)_2^{k,s}(c)$ for $k\geq 1,s\geq 2,c\in\{0,1\},k+s\geq 4.$

\item
If $A$ has three simple modules, then $A$ is stably equivalent of Morita type to
one and only one algebra of the type $SD(3\mathcal{K})^{a,b,c})$ for
$a\geq b\geq c\geq 1$.
\end{enumerate}

\item Let $A$ be an algebra of quaternion type.

\begin{enumerate}
\item If $A$ has one simple modules, then $A$ is stably equivalent
of Morita type to one of the algebras $Q(1{\mathcal A})_1^k$ for
$k\geq 2$ or $Q(1{\mathcal A})_2^k(c,d)$ for $k\geq 2, (c,d)\neq
(0,0)$ if characteristic if the $K$ is $2$. Different parameters
$k$ yield algebras in different
 equivalence classes of Morita type.

\item If $A$ has two simple modules then $A$ is stably equivalent of Morita
type to one of the algebras $Q(2{\mathcal B})_1^{k,s}(a,c)$
for $k\geq 1,s\geq 3,a\neq 0$.

\item If $A$ has three simple modules, then $A$ is stably equivalent of Morita type to
one of the algebras $Q(3{\mathcal K})^{a,b,c}$ for
$a\geq b\geq c\geq 1, b\geq 2, (a,b,c)\neq (2,2,1)$ or
$Q(3{\mathcal A})_1^{2,2}(d)$ for $d\in K\setminus\{0,1\}$.
Different parameters $a,b,c$ yield algebras in different stable equivalence
classes of Morita type.
\end{enumerate}

\end{itemize}
\end{Thm}

\begin{Proof}
The first point is Proposition~\ref{onlyinternal} and the second point is
Corollary~\ref{ARconjdihedral} and Corollary~\ref{ARconjfortametype}. The third point is
Proposition~\ref{dihedralonesimple}, Proposition~\ref{stableclassDiherdaltwosimples} and
Proposition~\ref{stableclassDiherdalthreesimples}.
The fourth point is Proposition~\ref{result} together with Section~\ref{semidihedralsection}
and the fifth point is Proposition~\ref{result}
together with Section~\ref{quaternionsection}.
\end{Proof}

\begin{Rem}\rm
For algebras of dihedral type, we proved in Section~\ref{dihedraltypesection}
that the
classification up to stable equivalences of Morita type coincide
with derived equivalence classification, up to a scalar problem in
$D(1\mathcal{A})_2^k(d)$. The only piece that is missing for a complete
classification is the question if
$D(1\mathcal{A})_2^k(0)$ is stably equivalent of
Morita type to $D(1\mathcal{A})_2^k(1)$.

Derived equivalent local algebras are Morita equivalent as is shown by
Roggenkamp and the second author (cf \cite{ZimmermannOrder}).
Observe that tame local symmetric algebras are classified in
\cite[Chapter III]{ErdmannLNM}. Actually, the classification coincides
with the algebras with one simple module we already dealt with in the text.
So, a complete classification of the algebras of dihedral type with one simple
module would give a
classification of tame local symmetric algebras.
\end{Rem}

\begin{Cor}
The Auslander Reiten conjecture holds for tame local symmetric algebras, i.e.
if $A$ is a tame local symmetric algebra and if $B$ is an algebra without simple direct
factor which is
stably equivalent of Morita type to $A$, then $B$ is local tame symmetric as well.
\end{Cor}

\begin{Proof}
By Liu \cite{Liu2008} $B$ is indecomposable since $A$ is indecomposable.
Erdmann classified tame local symmetric algebras
\cite[III.1 Theorem]{ErdmannLNM}. The classification coincides with
the list of local algebras of dihedral, semi-dihedral or quaternion type.
\end{Proof}

We cannot give any answer to the classification of algebras of
dihedral, semi-dihedral or quaternion type up to derived equivalence beyond
the information that is already known. Nevertheless, one more
statement for algebras of semi-dihedral type was obtained by Holm and the
second author.

\begin{Thm} (Holm and Zimmermann \cite{HolmZimmermann})\label{HZtame}
\begin{enumerate}
\item \label{thm1-intro-semidihedral} Let $F$ be an algebraically
closed field of characteristic 2. For any given integers $k,s\ge
1$, consider the algebras of semi-dihedral type $SD(2\mathcal{B})_1^{k,s}(c)$
for the scalars $c=0$ and $c=1$. Put
$B^{k,s}_c:=SD(2\mathcal{B})_1^{k,s}(c)$. Suppose that if $k=2$ then $s\ge
3$ is odd, and if $s=2$ then $k\ge 3$ is odd.
Then the factor rings $Z(B^{k,s}_0)/T_1(B^{k,s}_0)^{\perp}$ and
$Z(B^{k,s}_1)/T_1(B^{k,s}_1)^{\perp}$ are not isomorphic.

In particular, the algebras $SD(2\mathcal{B})_1^{k,s}(0)$ and
$SD(2\mathcal{B})_1^{k,s}(1)$ are not derived equivalent.

\item
\label{thm2-intro-semidihedral} Let $F$ be an algebraically closed
field of characteristic 2. For any given integers $k,s\ge 1$,
consider the algebras of semi-dihedral type $SD(2\mathcal{B})_2^{k,s}(c)$ for
the scalars $c=0$ and $c=1$. Put $C^{k,s}_c:=SD(2\mathcal{B})^{k,s}_2(c)$.
If the parameters $k$ and $s$ are both odd, then the factor rings
$Z(C^{k,s}_0)/T_1(C^{k,s}_0)^{\perp}$ and
$Z(C^{k,s}_1)/T_1(C^{k,s}_1)^{\perp}$ are not isomorphic. Hence
the algebras $SD(2\mathcal{B})_2^{k,s}(0)$ and $SD(2\mathcal{B})_2^{k,s}(1)$ have
different sequences of generalised Reynolds ideals.

In particular, for $k$ and $s$ odd, the algebras
$SD(2\mathcal{B})_2^{k,s}(0)$ and $SD(2\mathcal{B})_2^{k,s}(1)$ are not derived
equivalent.
\end{enumerate}
\end{Thm}

We get the following positive result.

\begin{Cor}\begin{enumerate}
\item Let $F$ be an algebraically closed field of characteristic
$2$. For any given integers $k,s\ge 1$, consider the algebras of
semi-dihedral type $SD(2\mathcal{B})_1^{k,s}(c)$ for the scalars $c=0$ and
$c=1$.  Suppose that if $k=2$
then $s\ge 3$ is odd, and if $s=2$ then $k\ge 3$ is odd. Then
the algebras $SD(2\mathcal{B})_1^{k,s}(0)$ and
$SD(2\mathcal{B})_1^{k,s}(1)$ are not stably equivalent of Morita type.

\item Let $F$ be an algebraically closed field of characteristic
$2$. For any given integers $k,s\ge 1$, consider the algebras of
semi-dihedral type $SD(2\mathcal{B})_2^{k,s}(c)$ for the scalars $c=0$ and
$c=1$. If the parameters $k$ and $s$ are both odd, then the
algebras $SD(2\mathcal{B})_2^{k,s}(0)$ and $SD(2\mathcal{B})_2^{k,s}(1)$ are not
stably equivalent of Morita type.
\end{enumerate}
\end{Cor}

\begin{Proof}
Since the quotients $Z^{st}(A):=Z(A)/Z^{pr}(A)$ and
$T_n^\perp(A)^{st}:=T_n(A)^\perp/Z^{pr}(A)$ are invariants under
stable equivalences of Morita type, so are the quotients
$Z^{st}(A)/T_n^\perp(A)^{st}=Z(A)/T_n^\perp(A)$.

Hence the parameters in the theorem yield not only algebras in
different derived equivalence classes, but also algebras in
different equivalence classes up to stable equivalences of Morita type.
\end{Proof}

\end{document}